\documentclass[11pt]{article}
\usepackage{amssymb}
\usepackage{amsfonts}
\usepackage{amsmath}
\usepackage{amsopn}

\def\sqr#1#2{{\vcenter{\vbox{\hrule height.#2pt
              \hbox{\vrule width.#2pt height#1pt \kern#1pt \vrule width.#2pt}
              \hrule height.#2pt}}}}
\def\signed #1{{\unskip\nobreak\hfil\penalty50
              \hskip2em\hbox{}\nobreak\hfil#1
              \parfillskip=0pt \finalhyphendemerits=0 \par}}
\def\endpf{\signed {$\sqr69$}}
\def\3n{\negthinspace \negthinspace \negthinspace }
\def\2n{\negthinspace \negthinspace }
\def\1n{\negthinspace }

\def\see{{\it see} }

\def\={\buildrel \triangle \over =}

\def\ds{\displaystyle}

\def\ns{\noalign{\ss}}

\def\a{\alpha}
\def\b{\beta}

\def\d{\delta}
\def\e{\varepsilon}

\def\l{\lambda}

\def\D{\Delta}

\def\cB{{\cal B}}

\def\cD{{\cal D}}

\def\cF{{\cal F}}

\def\cH{{\cal H}}

\def\cL{{\cal L}}

\def\cN{{\cal N}}

\def\cP{{\cal P}}
\def\cQ{{\cal Q}}
\def\cR{{\cal R}}

\def\cU{{\cal U}}

\def\ss{\smallskip}
\def\ms{\medskip}

\def\q{\quad}
\def\qq{\qquad}

\def\lan{\mathop{\langle}}
\def\ran{\mathop{\rangle}}

\def\pa{\partial}

\def\cd{\cdot}

\def\span{\hbox{\rm span$\,$}}

\def\deq{\mathop{\buildrel\D\over=}}

\def\({\Big (}
\def\){\Big )}
\def\[{\Big[}
\def\]{\Big]}
\def\bde{\begin{definition}}
\def\ede{\end{definition}}
\def\be{\begin{equation}}
\def\bel{\begin{equation}\label}
\def\ee{\end{equation}}
\def\bt{\begin{theorem}}
\def\et{\end{theorem}}
\def\bc{\begin{corollary}}
\def\ec{\end{corollary}}
\def\bl{\begin{lemma}}
\def\el{\end{lemma}}
\def\bp{\begin{proposition}}
\def\ep{\end{proposition}}
\def\bas{\begin{assumption}}
\def\eas{\end{assumption}}
\def\br{\begin{remark}}
\def\er{\end{remark}}
\def\ba{\begin{array}}
\def\ea{\end{array}}
\def\ed{\end{document}}

\def\square#1{\vbox{\hrule\hbox{\vrule height#1%
     \kern#1\vrule}\hrule}}
\def\rectangle#1#2{\vbox{\hrule\hbox{\vrule height#1%
     \kern#2\vrule}\hrule}}

\font\tenbb=msbm10 \font\sevenbb=msbm7 \font\fivebb=msbm5

\newfam\bbfam
\scriptscriptfont\bbfam=\fivebb \textfont\bbfam=\tenbb
\scriptfont\bbfam=\sevenbb

\newtheorem{lemma}{Lemma}[section]
\newtheorem{remark}{Remark}[section]

\newtheorem{theorem}{Theorem}[section]
\newtheorem{corollary}{Corollary}[section]

\newtheorem{definition}{Definition}[section]
\newtheorem{proposition}{Proposition}[section]
\newtheorem{assumption}{Assumption}[section]

\makeatletter
   \renewcommand{\theequation}{%
            \thesection.\arabic{equation}}
   \@addtoreset{equation}{section}
\makeatother

\begin{document}

\title{\bf  Equivalent Conditions on Periodic Feedback Stabilization  for Linear Periodic Evolution Equations}
\date{}
\author{Gengsheng Wang\thanks{
 School of Mathematics and Statistics, Wuhan University, Wuhan,
430072, China. (wanggs62@yeah.net) The author was partially
supported by the National Natural Science Foundation of China under grant
11161130003.} \q \q Yashan Xu\thanks{School of
Mathematical Sciences, Fudan University, KLMNS, Shanghai 200433,
China. (yashanxu@fudan.edu.cn) }}

\maketitle

\begin{abstract}
This paper studies the periodic feedback stabilization for a class of linear $T$-periodic evolution equations.
 Several equivalent conditions on the linear periodic feedback stabilization  are obtained.  These conditions are related with the following subjects: the  attainable subspace of the controlled evolution equation under consideration;
 the unstable subspace (of the evolution equation with the null control) provided by the Kato projection; the Poincar$\acute{e}$ map  associated with the  evolution equation with the null control; and two unique continuation properties for the dual equations  on different time horizons $[0,T]$ and $[0,n_0T]$ (where $n_0$ is the sum of algebraic multiplicities of  distinct unstable eigenvalues of the Poincar$\acute{e}$ map).  It is also proved that a $T$-periodic controlled evolution equation is linear $T$-periodic feedback sabilizable if and only if it is linear $T$-periodic feedback sabilizable with respect to a finite dimensional subspace.
Some applications to heat equations with time-periodic potentials are presented.
\end{abstract}

\ms


\ms


\rm

\section{Introduction}\label{1}

\subsection{The problem and the motivation}

Consider  the following controlled evolution equation:
\begin{equation}\label{state}
y'(t) + Ay(t) + B(t)y(t) = D(t)u(t)\;\;\mbox{in}\;\; \mathbb{R}^+\triangleq [0,\infty).
\end{equation}
Here and throughout this paper, we make the following assumptions.
\vskip 5pt

\noindent $(\cH_1)$  The operator $(-A)$, with its domain $\cD(-A)$,  generates a  $C_0$ compact semigroup
$\{S(t)\}_{t\ge0}$ in a real Hilbert space $H$ (identified with its dual) with its norm and inner product denoted by $\|\cdot\|$  and  $\lan\cd,\cd\ran$, respectively.

\vskip 5pt
\noindent $(\cH_2)$ The operator-valued  function $B(\cd)\in L^1_{loc}(\mathbb{R}^+;\cL(H))$ is $T$-periodic, i.e., $B(t+T)=B(t)$ for a.e. $t\in \mathbb{R}^+$, where $T>0$ and $\cL(H)$ denotes  the space of all linear
bounded operators on $H$.

\vskip 5pt

\noindent $(\cH_3)$ The operator-valued  function $D(\cd)\in L^\infty(\mathbb{R}^+;\cL(U,H))$ is $T$-periodic. Here $U$ is also a real Hilbert space (identified with its dual) with its norm and inner product denoted by $\|\cdot\|_U$ and  $\lan\cd,\cd\ran_U$, respectively; and  $\cL(U, H)$ stands for the space of all linear
bounded operators from $U$ to $H$. Controls $u(\cdot)$ are taken from the space $L^2(\mathbb{R}^+; U)$.

\vskip 5pt

\noindent
For each $h\in H$, $s\geq 0$ and $u(\cdot)\in L^2(\mathbb{R}^+ ; U)$, Equation (\ref{state}) (over $[s,\infty)$) with the initial condition that $y(s)=h$  has a  unique mild solution  $y(\cdot; s,h,u)\in C([s,\infty); H)$. (See, for instance, Proposition 5.3 on Page 66 in \cite{Li}.)
The following definitions about the periodic feedback stabilization will be used throughout this paper:

 \begin{itemize}

\item    Equation (\ref{state}) is said to be linear periodic feedback stablizable (LPFS, for short) if there is a $T$-periodic
$K(\cdot)\in L^\infty\left(\mathbb{R}^+;\cL(H, U)\right)$
 such that the feedback equation
 \be\label{fstate}
y'(t) + Ay(t) + B(t)y(t) = D(t)K(t)y(t)\quad\;\;\mbox{in}\;\;\mathbb{R}^+
\ee
 is exponentially stable, i.e., there are two positive constants $M$ and $\delta$ such that for each  $h\in H$, the solution $y_K(\cdot;0,h)$ to the equation (\ref{fstate}) with the initial condition that
    $y(0)=h$
satisfies that $\|y_K(t; 0,h)\|\leq M e^{-\d  t}\|h\|$ for all $t\geq 0$.
  Any such a $K(\cdot)$ is called an LPFS law for Equation (\ref{state}).

\item     Equation (\ref{state}) is said to be LPFS with respect to a subspace $Z$ of  $U$ if there is a $T$-periodic
$K(\cdot)\in L^\infty\left(\mathbb{R}^+;\cL(H, Z)\right)$ such that the equation (\ref{fstate}) is exponentially stable. Any such a $K(\cdot)$ is called an LPFS law for Equation  (\ref{state}) with respect to  $Z$.

 \end{itemize}
Let
\begin{equation}\label{space}
\cU^{FS}\triangleq \big\{ Z \; \big |\;  Z\;\;\mbox{is a subspace of}\;\; U\;\;\mbox{s.t. Equation (\ref{state}) is LPFS w.r.t.}\;Z \big\}.
\end{equation}

 In this paper, we provide three criteria for judging whether a subspace $Z$ belongs to $\cU^{FS}$. We also show that if $U\in \cU^{FS}$, then there is a  finite dimensional subspace $Z$ in $\cU^{FS}$.
 The aforementioned three criteria are related  with  the following subjects: the  attainable subspace of (\ref{state});
 the unstable subspace (of  (\ref{state}) with the null control) provided by the Kato projection; the Poincar$\acute{e}$ map  associated to (\ref{state}) with the null control; and two unique continuation properties for the dual equations of  (\ref{state}) (with the null control)  on different time horizons $[0,T]$ and $[0,n_0T]$ (where $n_0$ is the sum of algebraic multiplicities of  distinct unstable eigenvalues of the Poincar$\acute{e}$ map).
   Among three criteria, the  most important one  is a geometric condition  connecting  the attainable set   with the unstable subspace of the system (\ref{state});     while the other two  are analytic conditions related with the unique continuation of the dual equations of (\ref{state}) over different time horizons and with initial data in different finite dimensional subspaces
     of $H$.

The motivation for this work is as follows. First,  the equation (\ref{state})  with the null control is exponentially stable if and only if the spectrum of    the  Poincar$\acute{e}$ map associated with the system is contained in the unit open ball of the complex plane. (This can be proved by the exactly  same way to show  Corollary 7.2.4 on page 200, \cite{Henry}.)  Thus, it is a  natural  problem to explore equivalent conditions on the periodic stabilization for a linear periodic controlled evolution system. Second, there are two important kinds of solutions for evolution equations:  equilibrium and periodic solutions.
  The stabilization for equilibrium  solutions of time-invariant  systems has been extensively studied (see for instance \cite{Barbu2}, \cite{Coron}, \cite{Pritchard}, \cite{Russell} and the references therein). However,
 the understanding  on the periodic stabilization of periodic solutions for time-varying  evolution systems is quite limited. (See \cite{Barbu}, \cite{Lunardi}, \cite{Phung0} and \cite{Prato}. Here, we would like to mention
  \cite{Barbu0} which establishes a feedback law  stabilizing a smooth non-stationary  solutions, for instance, around a periodic trajectory, for Navier-Stokes equations.)
  Finally,  when the system (\ref{state}) is LPFS, it should be important and interesting to answer  if there is a finite dimensional subspace $Z$ of $U$ such that (\ref{state}) is LPFS w.r.t. $Z$, from perspectives of both mathematics and applied sciences.

 \subsection{Main results}

 Before stating our main  results, we give  some preliminaries  in order:

\noindent $(I)$  {\it Notations} $\;\;$
We will  use $\|\cdot\|$ to denote the usual norm of $\cL(H)$ when there is no risk of causing  any confusion.
Given  $L\in \cL(X,Y)$ (where $X$ and $Y$ are two Hilbert spaces),  we write $L^*$ for its adjoint operator. For $L\in\cL(X)\triangleq\cL(X,X)$, we denote by $\sigma(L)$  the spectrum of $L$. When $X$ is a real Hilbert space and $L\in\cL(H)$, we  denote  by $X^C$ and
 $L^C$  their complexification, respectively, i.e., $X^C=X+iX$ and $L^C(\a+i\b)=L\a+iL\b$ for any $\a, \b\in X$, where $i$ is the imaginary unit.
We write $\mathbb{B}$ for the open unit ball in $\mathbb{C}^1$ and $\mathbb{B}(0,\delta)$  for the open ball in $\mathbb{C}^1$, centered at the origin and of radius $\delta>0$.
Denote by $\partial\mathbb{B}(0,\delta)$ the boundary of $\mathbb{B}(0,\delta)$.

\noindent{(II)} {\it The Poincar$\acute{e}$ map}  $\;\;$
Let $\bigr\{\Phi(t, s)\bigl\}_{0\le s\le t<+\infty}$ be the evolution system generated by
$(-A-B(\cd))$. It follows from  Lemma 5.6   in  \cite{Li} (see Page 68, \cite{Li}) that $\Phi(t,s)$ is strongly continuous over $\{(t,s)\in \mathbb{R}^+\times \mathbb{R}^+ \big | 0\leq s\leq t<\infty\}$,
and that
\begin{equation}\label{huang1.4}
\Phi(t,s)h=S(t-s)h+\ds\int_s^tS(t-r)B(r)\Phi(r,s)hdr,\;\;\mbox{when}\;\; 0\leq s\leq t<\infty\;\;\mbox{and}\;\; h\in H.
\end{equation}
By Proposition 5.7  on Page 69  in \cite{Li}, for each $s\geq 0$,   $h\in H$  and $u(\cdot)\in L^2(\mathbb{R}^+;U)$, it holds that
\begin{equation}\label{solution1}
y(t;s,h, u)=\Phi(t,s)h+\int^t_s\Phi(t,r) D(r)u(r)dr, \;\;  s\in [t, \infty).
\end{equation}
By the $T$-periodicity of $B(\cd)$ and (\ref{huang1.4}), one can easily check that  $\Phi(\cd, \cd)$ is  $T$-periodic, i.e.,
\begin{equation}\label{wgs1.5}
\Phi(t+T, s+T)=\Phi(t,s)\;\; \mbox{ for all }\;0\le s\le t<\infty.
\end{equation}
Now, we introduce the following Poincar$\acute{e}$ map  (\see  Page 197, \cite{Henry}):
\begin{equation}\label{poincare}
\cP(t) \triangleq\Phi(t+T,t),\;\;t\in\mathbb{R}^+.
\end{equation}
It is proved that (see  Lemma \ref{lemma2.1})
\begin{equation}\label{WGS1.8}
\sigma(\cP(t)^C)\setminus\{0\}=\{\l_j\}^\infty_{j=1}\;\;\mbox{for each}\;\; t\geq 0,
\end{equation}
where $\l_j$, $j=1,2,\dots$, are all distinct non-zero eigenvalues of the compact operator $\cP(0)^C$ such that $\lim_{j\rightarrow\infty}|\l_j|=0$. Thus, there is a unique $n\in \mathbb{N}$ such that
 \begin{equation}\label{wgs1.9}
 |\l_j|\ge1,\q j\in\{1,2,\cdots, n\}\q~{\mbox{and}}~\q |\l_j|<1,\q j\in\{n+1,n+2,\cdots \}.
 \end{equation}
Set
 \begin{equation}\label{decay}
 \bar \delta\triangleq \max\{|\l_j|,\,  j>n\}<1.
 \end{equation}
 Let $l_j$ be the algebraic multiplicity of $\l_j$ for each $j\in \mathbb{N}$, and write
 \begin{equation}\label{number}
n_0\triangleq l_1+ \cdots +l_n.
 \end{equation}
\noindent $(III)$ {\it The Kato projection} $\;\;$  Arbitrarily fix a $\delta\in (\bar \delta, 1)$, where $\bar \delta$ is given by (\ref{decay}).   Let
$\Gamma$ be the circle $\pa \mathbb{B}\left(0, \d\right)$ with the clockwise direction in $\mathbb{C}^1$.
We introduce the  Kato projections (see \cite{Kato}):
\begin{equation}\label{wgs1.11}
\hat P(t) =\ds\frac{1}{ 2\pi i}\int_\Gamma (\l I-\cP(t)^C)^{-1}d\l,\;\; t\geq 0.
\end{equation}
It is proved that (see Lemma~\ref{lemma2.2}) for each $t\geq 0$,
the operator $P(t)$, defined by
\begin{equation}\label{P}
P(t)\triangleq\hat P(t)\bigm|_H\;\;(\mbox{the restriction of}\;\hat P(t)\;\mbox{on}\; H),
\end{equation}
is a projection on $H$; $H= H_1(t)\bigoplus H_2(t)\triangleq P(t)H\bigoplus (I-P(t))H$ for each $ t\geq 0$;
  both $H_1(t)$ and $H_2(t)$ are invariant subspaces of $\cP(t)$; $\sigma(\cP(t)^C|_{H_1(t)^C}) =\{\l_j \}^n_{j=1}$, $\sigma(\cP(t)^C|_{H_2(t)^C})\setminus\{0\} = \{\l_j \}^\infty_{j=n+1}$; and $\mbox{dim} H_1(t)=n_0$.
It is also shown  that (see Lemma~\ref{lemma2.2})  $P(\cdot)$ is $T$-periodic.
 We simply write
\begin{equation}\label{wgs1.15}
H_1\triangleq H_1(0),\; H_2\triangleq H_2(0),\; P\triangleq P(0)\;\;\mbox{and}\;\;\cP\triangleq\cP(0).
\end{equation}
 The subspaces $H_1$ and $H_2$ are respectively called the unstable subspace and  the stable subspace of Equation (\ref{state}) with the null control.  Each  eigenvalue in $\{\l_j\}_{j=1}^n$ (or in $\{\l_j\}_{j=n+1}^\infty$) is called an unstable (or stable) eigenvalue of $\cP^C$. Each eigenfunction of $\cP^C$ corresponding to an unstable (or stable) eigenvalue is called an unstable (or stable) eigenfunction of $\cP^C$.

\noindent $(IV)$ {\it Attainable subspaces} $\;\;$  For each subspace $Z\subset U$, we let
\begin{equation}\label{csubspace}
V^Z_{k}\deq \Big\{\ds\int^{kT}_0\Phi(kT,s)D(s)  u(s) ds~\Bigm| u(\cdot)\in L^2(\mathbb{R}^+;Z)\Big\}\;\;{\hbox { for all }}\;\; k\in\mathbb{N}.
\end{equation}
The space $V^Z_{k}$ is called the attainable subspace of Equation (\ref{state}) (over $(0,kT)$) w.r.t. $Z$.
Let
\begin{equation}\label{2-14}
 \hat V^Z _{k}=P V^Z _{k},
\qq k\in\mathbb{N},
\end{equation}
where $P$ is given by (\ref{wgs1.15}).

\bigskip

 Now the main results of this paper are presented by the following two theorems:

\bt\label{theorem1} Let $P$, $\cP$ and $H_j$ with $j=1,2$ be given by (\ref{wgs1.15}). Let  $n_0$  be given by   (\ref{number}). Then, for each subspace $Z\subseteq U$,
   the following statements are equivalent:

  \noindent $(a)$  Equation (\ref{state}) is LPFS with respect to $Z$, i.e., $Z\in\cU^{FS}$.

  \noindent $(b)$  The subspace $Z$  satisfies
\begin{equation}\label{3-14}
\hat V^Z_{n_0}=H_1,\;\;\mbox{where}\;\;\hat V^Z_{n_0}\;\;\mbox{is given by}\;\; (\ref{2-14}).
\end{equation}

  \noindent  $(c)$ The subspace $Z$ satisfies
\begin{equation}\label{a.1.1}
 \xi\in P^*H_1\;\;\mbox{and}\;\; \left(D(\cdot)\bigl|_Z\right)^*\Phi(n_0T, \cdot)^*\xi=0\; \mbox{over}\;(0,n_0T)\Rightarrow \xi=0.
\end{equation}

\noindent  $(d)$ The subspace $Z$ satisfies
\begin{equation}\label{a.1.2}
\mu\notin\mathbb{B},\;\xi\in H^C,\;  \big(\mu I-\cP^{*C}\big)\xi=0,\;\mbox{ and }\;  \big(D(\cd)\bigl|_Z\big)^{*C}\Phi(T,\cd)^{*C}\xi=0\;\mbox{over}\; (0,T)\;\Rightarrow \xi=0.
\end{equation}

  \et

\bt\label{theorem2} Equation (\ref{state})  is LPFS if and only if it is LPFS with respect to a finite dimensional subspace $Z$ of $U$.
\et

 \medskip

 It is worthwhile to make the following remarks:

  \begin{itemize}

\item The key to show the above two theorems is to build up the equivalence $(a)\Leftrightarrow (b)$ in Theorem~\ref{theorem1}.

 \item The functions $\Phi(n_0T,\cdot)^*\xi$ with $\xi\in H$ and $\Phi(T,\cdot)^{*C}\xi$ with $\xi\in H^C$ are respectively  the solutions to the following dual equations:
 \begin{equation*}\label{yuhuang1.25}
 \psi_t(t)-A^*\psi(t)-B(t)^*\psi(t)=0\;\;\mbox{for a.e.}\; t\in (0,n_0T),\;\;\psi(n_0T)=\xi
 \end{equation*}
 and
 \begin{equation*}\label{yuhuang1.26}
\psi_t(t)-A^{*C}\psi(t)-B(t)^{*C}\psi(t)=0\;\;\mbox{for a.e.}\; t\in (0, T),\;\;\psi(T)=\xi.
 \end{equation*}
 Thus, the condition $(c)$ in Theorem~\ref{theorem1} presents a unique continuation property for solutions of the first dual equation   with initial data in $P^*H_1$;  while the condition $(d)$  in Theorem~\ref{theorem1} presents a unique continuation property for solutions
 of the second equation where initial data are unstable eigenfunctions of $\cP^{*C}$.

\item There have been studies, in the past,  on  equivalence conditions of periodic feedback stabilization for linear periodic evolution systems.
   In \cite{Lunardi} and \cite{Prato}, the authors established  an equivalent condition on stabilizability for linear time-periodic parabolic equations with open-loop controls.
 Their
  equivalence (see Theorem 3.1 in \cite{Lunardi} and Proposition 3.1 in \cite{Prato}) can be stated, under our framework,  as follows: the condition $(d)$ (in our Theorem~\ref{theorem1} where $Z=U$) is equivalent to
 the statement that for any $h\in H$, there is a control $u^h(\cdot)\in C(\mathbb{R}^+; U)$, with $\sup\limits_{t\in\mathbb{R}^+}\|e^{\bar\delta t}u^h(t)\|_U$ bounded (where $\bar\delta$ is given by (\ref{decay})), such that the solution $y(\cdot; 0,h, u^h)$ is stable. Meanwhile, it was pointed out in \cite{Lunardi} (see the paragraph before the last one in Section 1 in \cite{Lunardi}) that when open-looped stabilization controls exist, one can construct a periodic feedback stabilization law through using a method provided  in \cite{Prato1}.
 From this point of view, the equivalence $(a)\Leftrightarrow (d)$  in  Theorem~\ref{theorem1} is not new,  though our way to approach the equivalence differs  from those in \cite{Lunardi} and \cite{Prato} and our method to construct the stabilization feedback law is different from that in \cite{Prato1}.

   \item To the best of our knowledge,  both Theorem~\ref{theorem2} and the equivalences:  $(a)\Leftrightarrow (b)$ and $(a)\Leftrightarrow (c)$ in   Theorem~\ref{theorem1} appear to be new.
 It is worth mentioning that the  equivalence  $(a)\Leftrightarrow (b)$ in  Theorem~\ref{theorem1} is an extension of a result in our previous paper  \cite{Yashan} which studies the stabilization of  finite-dimensional periodic systems.

\item A byproduct of this study (see Proposition~\ref{cor3.1}, and Remark~\ref{remark3.1}) shows that when  both $B(\cdot)$ and $D(\cdot)$ are time-invariant, linear time-period
functions $K(\cdot)$ will not aid the linear stabilization of Equation (\ref{state}), i.e., Equation (\ref{state}) is linear $\hat T$-periodic feedback stabilizable for some $\hat T>0$ if and only if Equation (\ref{state}) is linear time-invariant feedback stabilizable. On the other hand, when Equation (\ref{state}) is periodic time-varying, linear time-periodic $K(\cdot)$ do aid in the linear stabilization of  this equation.

 \end{itemize}

  The rest of this paper is organized as follows.  Section 2 provides some properties on Poincar$\acute{e}$ map,  Kato projection and  attainable subspaces.
 Section 3 studies the multi-periodic feedback stabilization.  Section 4  proves  Theorem {\ref{theorem1} and {\ref{theorem2}. Section 5  presents some applications of the main theorems to internally controlled heat equations with time-periodic potentials.

 \section{  Poincar$\acute{e}$ map,  Kato projection and  attainable subspaces}

In this section, we will first present  three lemmas (Lemma 2.2, 2.3, 2.4) which show certain properties on  the Poincar$\acute{e}$ map and the Kato projection. (These lemmas  are slightly  different versions of the existing results. For the sake of the completeness of the paper, we provide their detailed proof  in Appendix.) Then, we show certain properties on attainable subspaces $V^Z_k$ defined by (\ref{csubspace}).

 The first one is another version of Lemma 7.2.2 in \cite{Henry}  where  $B(\cdot)$ is assumed to be $T$-periodic and  H$\ddot{o}$lder continuous (see Page 197, \cite{Henry}); while in our case, $B(\cdot)\in L^1(\mathbb{R}^+;\cL(H))$ is $T$-periodic.
\begin{lemma}\label{lemma2.1}
 Let  $\cP(\cdot)$ be defined by (\ref{poincare}).
Then $\sigma(\cP(t)^C)\setminus\{0\}$ is independent of $t\in \mathbb{R}^+$. Moreover, $\sigma(\cP(t)^C)\setminus\{0\}$ consists entirely of distinct eigenvalues $\{\l_j\}_{j=1}^\infty$  (of $\cP(0)$) with $\lim_{j\rightarrow\infty}|\l_j|=0$.
\end{lemma}
The second one is another version of Theorem  7.2.3 in \cite{Henry} where complex case is studied (see   Page 198, \cite{Henry}); while in following lemma, we consider the real case.
\begin{lemma}\label{lemma2.2}
Let   $\cP(\cdot)$  and $P(\cdot)$ be defined by (\ref{poincare}) and (\ref{P}) respectively.
Then  each  $P(t)$ (with $t\geq 0$) is a projection  on $H$
such that
\begin{equation}\label{2-2}
H=H_1(t)\bigoplus H_2(t),
\end{equation}
where
\begin{equation}\label{2-2-1}
H_1(t)\deq P(t)H\;\;\mbox{and}\;\;  H_2(t)\deq(I-P(t))H.
\end{equation}
Moreover, $P(\cdot)$, $H_1(\cdot)$ and $H_2(\cdot)$ have the following properties:

\noindent $(a)$ $P(\cdot)$, $H_1(\cdot)$ and $H_2(\cdot)$ are $T$-periodic.

\noindent $(b)$ For each $t\geq 0$, both $H_1(t)$ and $H_2(t)$ are invariant subspaces of $\cP(t)$.

\noindent $(c)$ If $\{\lambda_j\}_{j=1}^\infty$, $n$ and $n_0$ are given by (\ref{WGS1.8}), (\ref{wgs1.9}) and (\ref{number}), then
\begin{equation}\label{{2-3}}
\sigma(\cP(t)^C|_{H_1(t)^C}) =\{\l_j \}^n_{j=1},\qq \sigma(\cP(t)^C|_{H_2(t)^C})\setminus\{0\} = \{\l_j \}^\infty_{j=n+1};
\end{equation}
\begin{equation}\label{wgs1.14}
\mbox{dim} H_1(t)=n_0.
\end{equation}

\noindent $(d)$  When $0\le s\le t < \infty$, $\Phi(t, s) \in\cL(H_j(s), H_j(t))$ with $j=1,2$.

\noindent $(e)$ It holds that
\begin{equation}\label{2-4}
\Phi(t, s)P(s) = P(t)\Phi(t, s),\;\; \mbox{when}\;\; 0\le s\le t.    \end{equation}

\noindent $(f)$ Let $\bar \rho\triangleq (-{\ln \bar\d})/{T}>0$ with $\bar\delta$ given by (\ref{decay}).
For any $\rho\in(0,\bar\rho)$, there is a positive constant $ C_\rho$ such that
\begin{equation}\label{2-5}
\|\Phi(t, s)h_2\|\le  C_\rho e^{-\rho(t-s)}\|h_2\|, \mbox{ when}\;\;  0\le s\le t<\infty\;\;\mbox{and}\;\; h_2\in  H_2(s).
\end{equation}
\end{lemma}
The third one is essentially another version of Theorem 6.22 in \cite{Kato} (see Page 184, \cite{Kato}). To state it, we recall that  $\cP^{*}$ and $P^*$ are the adjoint operators of $\cP$ and $P$.
It is clear that
 \begin{equation}\label{a.1.41}
 \sigma\big(\cP^{*C}\big)=\overline{\sigma\big(\cP^C\big)}.
 \end{equation}
Since $\sigma\big(\cP^C\big)\setminus\{0\}=\{\l_j\}^\infty_{j=1}$ (see (\ref{WGS1.8})), it holds that $\sigma\big(\cP^{*C}\big)\setminus\{0\}=\{\bar\l_j\}^\infty_{j=1}$. Write
$\bar l_j$ for the algebraic multiplicity of $\bar\l_j$ w.r.t. $\cP^{*C}$. It is clear that
\begin{equation}\label{a.1.42}
\bar l_j=l_j\;\;\mbox{for all}\;\; j,\;\;\mbox{and}\;\; \bar l_1+\cdots +\bar l_n=n_0,
\end{equation}
where $l_j$ is the  algebraic multiplicity of $\l_j$ w.r.t. $\cP^C$;  $n$ and $n_0$ are given by (\ref{wgs1.9}) and (\ref{number}) respectively.
 Let $\Gamma$ be the circle used to defined $\hat P\triangleq\hat P(0)$ (see (\ref{wgs1.11})). Define the Kato projection with respect to $\cP^{*C}$ as follows:
 \begin{equation}\label{a.1.43}
\hat {\tilde P} =\ds\frac{1}{ 2\pi i}\int_\Gamma (\l I-\cP^{*C})^{-1}d\l.
\end{equation}
From Theorem 6.17 on Page 178 in \cite{Kato}, it follows that
\begin{equation}\label{YUHEHHEN2.22}
H^C=\hat {\tilde H}_1\bigoplus \hat {\tilde H}_2\triangleq \hat {\tilde P} H^C\bigoplus (I-\hat {\tilde P})H^C
\end{equation}
and
\begin{equation}\label{YUHEHHEN2.23}
\mbox{both}\;\; \hat {\tilde H}_1\;\mbox{and}\; \hat {\tilde H}_2\;\;\mbox{are invariant w.r.t.}\;\;\cP^{*C}.
\end{equation}

\begin{lemma}\label{HUANGlemma2.3}
Let $\hat {\tilde P}$ be defined by (\ref{a.1.43}). Then, $\tilde P\triangleq \hat {\tilde P}\big|_H$ is a projection on $H$;
$H=\tilde H_1\bigoplus\tilde H_2$, where $\tilde H_1\triangleq\tilde PH$ and $\tilde H_2\triangleq(I-\tilde P)H$;
$\cP^*\tilde H_1\subseteq \tilde H_1$;
$\sigma\big(\cP^{*C}|_{\tilde H_1^C}\big)=\{\bar\l_j\}^n_{j=1}$ and $\sigma\big(\cP^{*C}|_{\tilde H_2^C}\big)\subseteq\mathbb{B}$;
and $\mbox{dim}\tilde H_1=n_0$. It further holds that
\begin{equation}\label{a.1.49}
\tilde P=P^*;
\end{equation}
\begin{equation}\label{a.1.50}
\tilde H_1= P^*H=P^*H_1,\;\;\mbox{where}\;\; H_1\;\;\mbox{is given by}\;\;(\ref{wgs1.15});
\end{equation}
\begin{equation}\label{WANGJIAXU3.31}
\xi\in \tilde H_1^C,\;\;\mbox{when}\;\; \mu\in \sigma\big(\cP^{*C}\big)\setminus\mathbb{B}\;\;\mbox{and}\;\; (\mu I-\cP^{*C})\xi=0.
\end{equation}

\end{lemma}

\vskip 5pt

Next, we will introduce certain properties on attainable subspaces $V^Z _{k}$, $k\in\mathbb{N}$. They will play important roles in the proof of our main theorems.
We start with recalling  (\ref{2-2}) and (\ref{csubspace}).
 Since $H_1$ is invariant w.r.t. $\cP$ (see Part $(b)$ of Lemma~\ref{lemma2.2}), we can define $\cP_1: H_1\rightarrow H_1$ by setting
\begin{equation}\label{2-13}
\cP_1\triangleq{\large{\mbox{$\cP$}}}\bigm|_{H_1}.
\end{equation}
By  (\ref{{2-3}}), it holds that
\begin{equation}\label{WGSwgs2.22}
\sigma (\cP_1)\bigcap \mathbb{B}=\emptyset.
\end{equation}

\begin{lemma}\label{lemma2.3} Let  $\cP_1$  and $n_0$ be given by (\ref{2-13}) and (\ref{number}), respectively. Suppose that $Z\subseteq U$ is a subspace with  $ V^Z _{k}$ and $ \hat V^Z _{k}$
given by (\ref{csubspace}) and (\ref{2-14}), respectively.
Then for each $k\in\mathbb{N}$,
\begin{equation}\label{2-16}
 V^Z _k= V^Z _1+\cP V^Z _1+\cdots+\cP^{k-1} V^Z _1;\qq \hat  V^Z _k=\hat  V^Z _1+\cP_1\hat V^Z _1+\cdots+\cP_1^{k-1}\hat V^Z _1.
\end{equation}
Furthermore,  $\cP_1$ is invertible and it holds that
\begin{equation}\label{2-17}
\hat  V^Z =\hat  V^Z _{n_0};\;\;\;\; \cP_1\hat  V^Z =\hat V^Z =\cP_1^{-1}\hat  V^Z ,
\end{equation}
where
\begin{equation}\label{2-15}
  \hat  V^Z\triangleq\bigcup_{k=1}^\infty\hat V^Z_k.
\end{equation}
\end{lemma}

\noindent {\it Proof.} We begin with proving the first equality in (\ref{2-16}) by the mathematical induction.
Clearly, it stands when $k=1$.  Assume that it holds in the case when
$k=k_0$ for some $k_0\ge1$, i.e.,
\begin{equation}\label{wgs2.26}
V^Z_{k_0}=V^Z_1+\cP V^Z_1+\cdots +\cP^{k_0-1}V^Z_1.
\end{equation}
Because of (\ref{wgs1.5}) and (\ref{wgs1.15}), we have that $\Phi((k_0+1)T,T)=\Phi(T,0)^{k_0}=\cP^{k_0}$. This, along with (\ref{csubspace}), the $T$-periodicity of $D(\cdot)$ and (\ref{wgs2.26}), indicates that
$$\ba{rl} V^Z_{k_0+1}&=\Big\{\ds \int^{(k_0+1)T}_{0}\Phi((k_0+1)T,s)D(s)  u(s)ds \bigm|u(\cdot)\in L^2(\mathbb{R}^+; Z)\Big\}\\
&=\Big\{\Phi((k_0+1)T,T)\ds \int^{T}_{0}\Phi(T,s)D(s) u(s)ds\bigm|u(\cdot)\in L^2(\mathbb{R}^+; Z)\Big\}\\
&\q+\Big\{\ds\int^{ k_0T}_0\Phi((k_0+1)T,s+T)D(s+T) u(s+T)ds\bigm|u(\cdot)\in L^2(\mathbb{R}^+; Z)\Big\}\\
&=\cP^{k_0} V^Z _1+
\Big\{\ds\int_{0}^{k_0T}\Phi(k_0T,s)D(s)  u(s+T)ds\bigm|u(\cdot)\in L^2(\mathbb{R}^+; Z)\Big\}\\
\ns&=\cP^{k_0} V^Z _1+ V^Z_{k_0}= V^Z _1+\cP V^Z _1+\cdots+\cP^{k_0} V^Z _1.\ea
$$
which leads to the first equality in (\ref{2-16}).

We next show the second equality in (\ref{2-16}). By (\ref{wgs1.15}) and (\ref{2-4}) with $t=T$ and $s=0$, we have
\begin{equation}\label{wgs2.27}
\cP P=P\cP.
\end{equation}
Since $ P$ is a projection from $H$ onto $ H_1$ (see Lemma~\ref{lemma2.2}),  the second equality in (\ref{2-16}) follows from the first one in (\ref{2-16}), (\ref{2-14}) and (\ref{wgs2.27}).

Then we show the first equality in (\ref{2-17}). It follows respectively from  (\ref{2-15}) and (\ref{2-16}) that
\begin{equation}\label{wgs2.29}
\hat V^Z _{n_0}\subseteq \hat V^Z\;\;\mbox{and}\;\;\hat V^Z_k\subseteq \hat V^Z_{n_0},\;\;\mbox{when}\;\; k\leq n_0.
\end{equation}
Since $\mbox{dim}H_1=n_0$ (see  (\ref{wgs1.14})) and $\cP_1: H_1\rightarrow H_1$ (see (\ref{2-13})), according to
  the Hamilton-Cayley theorem, each
 $\cP_1^{j}$ with $ j\ge n_0$ is a linear combination of  $\big\{~I,~\cP_1^{1},~\cP_1^{2},\cdots,~\cP_1^{(n_0-1)}\big\}$.
This, along with the second equality in (\ref{2-16}),  indicates that
\begin{equation}\label{wgs2.30-1}
\hat V^Z _{k}=\sum_{j=0}^{k-1}\cP_1^{j}(\hat V^Z_1)\subseteq\sum_{j=0}^{n_0-1} \cP_1^{j} (\hat V^Z_1)=\hat V^Z _{n_0},\;\;\mbox{when}\;\; k\ge n_0.
\end{equation}
Now the first equality in (\ref{2-17}) follows from  (\ref{wgs2.29})  and (\ref{wgs2.30-1}).

Finally, we show  the non-singularity of $\cP_1$ and the second equality  in (\ref{2-17}). By the first equality in (\ref{2-17}) and the  Hamilton-Cayley theorem, we see that
$$\cP_1 \hat V^Z  =\cP_1\hat V^Z _{n_0}=\cP_1\sum_{j=0}^{n_0-1}\cP_1^{j} (\hat V^Z_1)
=\sum_{j=1}^{n_0}\cP_1^{j} (\hat V^Z_1)\subseteq\sum_{j=0}^{n_0-1}\cP_1^{j} (\hat V^Z_1)=\hat V^Z _{n_0},$$
from which, it follows that
\begin{equation}\label{2-18}
\cP_1 \hat V^Z \subseteq \hat V^Z .
\end{equation}
Because $0\notin\sigma(\cP_1^C)$ (see (\ref{{2-3}}) as well as (\ref{wgs1.9})), and since the domain of $\cP_1^C$,   $ H_1^C$ is a  finite dimensional subspace, the operator $\cP_1^C$ is invertible. Hence, the operator  $\cP_1$ is also invertible. This implies that  $\mbox{dim} (\cP_1 \hat V^Z )= \mbox{dim} \hat V^Z $, which,  together with
(\ref{2-18}), yields that $\cP_1 \hat V^Z = \hat V^Z $.
This completes the proof.  \endpf

\begin{lemma}\label{lemma3.2}
Let   $n_0$ be given by (\ref{number}). Then, for each subspace  $Z$ of $U$,  there is a finite dimension subspace of $\hat Z$ of $Z$ such that
\begin{equation}\label{3-18}
\hat V^Z_{n_0}=\hat V^{\hat Z}_{n_0},
\end{equation}
where $\hat V^Z_{n_0}$ and $\hat V^{\hat Z}_{n_0}$ are defined by (\ref{2-14}).
\end{lemma}
\noindent{\it Proof.} Let $Z$ be a subspace of $U$. Since $\hat V^Z_{n_0}$ is a subspace of $H_1$ and $\mbox{dim}H_1=n_0<\infty$ (see (\ref{wgs1.14}), we can assume that
$\mbox{dim}\hat V^Z_{n_0}\triangleq m\le n_0$. Write $\{\xi_1,\dots,\xi_m\}$  for
 an orthonormal basis of $\hat V^Z_{n_0}$. By the definition of $\hat V^Z_{n_0}$ (see (\ref{2-14}), as well as (\ref{csubspace})), there are $u_j(\cdot)\in L^2(\mathbb{R}^+;Z)$, $j=1,\dots,m$,
 such that
\begin{equation}\label{wgs3.39}
\ds \ds\int^{n_0T}_0P\Phi(n_0T, s)D(s)u_j(s)ds=\xi_j\;\;\mbox{ for all }\;\;j=1,\dots,m.
\end{equation}
By the boundedness of $\Phi(n_0T,\cdot)$  and $D(\cd)$ over $[0,n_0T]$, there is a constant $C>0$ such that
\begin{equation}\label{wgs3.40}
\|P\Phi(n_0T,\cdot)D(\cdot)\|_{L^\infty(0,n_0T;\mathcal{L}(U;H))}\leq C.
\end{equation}
Let $\hat \varepsilon>0$ small enough to satisfy
 \begin{equation}\label{HUANG3.41}
 (1-\hat\e)^2-(2\hat\e+\hat\e^2)(m-1) >0.
 \end{equation}
By the definition of the Bochner integration (\see \cite{Bochner}),  there are simple functions
 \begin{equation}\label{WHUANG3.42}
  v_j(\cdot)=\sum\limits_{l=1}^{k_j}{\mbox{\Large $\chi$}}_{E_{j\hspace{0.5mm}l}}(\cdot)z_{j\hspace{0.5mm}l}\;\;\mbox{over}\;\;(0,n_0T),\;\;j=1,\dots,m,
  \end{equation}
with  $z_{jl}\in Z$, $E_{jl}$  measurable sets in $(0,n_0T)$ and   $\chi_{E_{jl}}$  the characteristic function of  ${E_{jl}}$,
 such that
 $$
 \ds\int_0^{n_0T}\|u_j(s)- v_j(s)\|_Uds\leq \hat \varepsilon/C.
 $$
 This, along with (\ref{wgs3.40}), yields that for each $j\in\{1,\dots,m\}$,
\begin{equation}\label{wgs3.41}
 \ds\Big\|\int^{n_0T}_0P\Phi(n_0T, s)D(s)u_j(s)ds-  \ds\int^{n_0T}_0P\Phi(n_0T, s)D(s) v_j(s)ds\Big\|\le\hat\e.
 \end{equation}
 Let
\begin{equation}\label{HUANGWANGHEN3.44}
\eta_j=\ds\int^{n_0T}_0P\Phi(n_0T, s)D(s) v_j(s)ds,\;\;j=1,\dots,m.
\end{equation}
By  (\ref{2-14}), (\ref{csubspace}), (\ref{WHUANG3.42}) and (\ref{HUANGWANGHEN3.44}), we see that  $\eta_j\in \hat V^Z_{n_0}$ for all $j=1,\dots,m$. Meanwhile, it follows from (\ref{wgs3.41}) that
\begin{equation}\label{3-18-1}
\|\eta_j-\xi_j\|\le\hat\e\;\;\mbox{ for all }\;\;j\in\{1,\dots, m\}.
\end{equation}
Now we claim that $\{\eta_1,\dots,\eta_m\}$ is a basis of $\hat V^Z_{n_0}$.
 In fact, since $\{\xi_1,\dots,\xi_m\}$ is orthonormal, it follows from (\ref{3-18-1}) that
$$
|\lan \eta_j,\eta_l\ran|=| \lan\xi_j+(\eta_j-\xi_j),\xi_l+(\eta_l-\xi_l)\ran| \le 2\hat\varepsilon+\hat\varepsilon^2\;\;\mbox{for all}\;\;j,l\in\{1,\dots,m\}\;\;\mbox{with}\;\;j\neq l.
$$
From this and  (\ref{3-18-1}),  one can directly check that
$$\ba{rl}
\ns &\ds \lan \sum^{m}_{j=1}\alpha_j\eta_j,\,\sum^m _{l=1}\alpha_l\eta_l \ran
=\ds \sum^{m}_{j=1}\alpha_j^2\|\eta_j\|^2+\sum^{m}_{j=1}\sum _{l\neq j}\alpha_l\alpha_j\lan\eta_j,\,\eta_l\ran\\
\ns\ge&\ds\sum^{m}_{j=1}\alpha_j^2\left(\|\xi_j\|-\|\xi_j-\eta_j\|\right)^2-\sum^{m}_{j=1}\sum _{l\neq j}|\alpha_l\alpha_j|\cdot|\lan\eta_j,\,\eta_l\ran|\\
 \ns\geq &\ds\left((1-\hat\e)^2-(2\hat\e+\hat\e^2)(m-1)\right)\sum^{m}_{j=1}\alpha_j^2,\;\;\mbox{when}\;\; \alpha_1,\dots,\alpha_m\in \mathbb{R}.
 \ea$$
This, along with (\ref{HUANG3.41}), indicates that $\alpha_1=\cdots=\alpha_m=0$ whenever $\sum^{m}_{j=1}\alpha_j\eta_j=0$. Namely, $\{\eta_1,\dots,\eta_m\}$ is linearly independent group in
the subspace $\hat V^Z_{n_0}$ which has the dimension $m$. Hence, $\{\eta_1,\dots,\eta_m\}$ is a basis of $\hat V^Z_{n_0}$.

Let
$$\hat Z=\span\left\{z_{11},\dots,z_{1k_1},z_{21},\dots,z_{2k_2},\dots,z_{m1},\dots,z_{mk_m}\right\},$$
where $z_{jl}$, $j=1,\dots,m$, $l=1,\dots, k_j$, are given by (\ref{WHUANG3.42}). Clearly, $\hat Z$ is a finite-dimensional subspace of $Z$ and all $v_j(\cdot)$, $j=1,\dots, m$, (given by (\ref{WHUANG3.42})) belong to $L^2(\mathbb{R}^+;\hat Z)$.
These, along with (\ref{HUANGWANGHEN3.44}), yield that all $\eta_j$, $j=1,\dots,m$, are in $ \hat V^{\hat Z}_{n_0}$.
Therefore, it holds that
$$
\hat V^{Z}_{n_0}\supseteq \hat V^{\hat Z}_{n_0}\supseteq\span\{\eta_1,\dots,\eta_m\}= \hat V^{ Z}_{n_0}.
$$
This leads to  (\ref{3-18}) and completes the proof.\endpf

\section{The multi-periodic feedback stabilization}

In this section, we will introduce three propositions. The first two propositions will be used in the proof of our main theorems. The last one is independent interesting.
 We begin with the following definitions:

 \begin{itemize}
\item  Equation (\ref{state}) is said to be linear multi-periodic feedback stabilizable (LMPFS, for short) if there is a $kT$-periodic
$K(\cdot)\in L^\infty\left(\mathbb{R}^+;\cL(H, U)\right)$ for some $k\in \mathbb{N}$ such that the equation (\ref{fstate}) is exponentially  stable. Any such a $K(\cdot)$ is called an LMPFS law for Equation (\ref{state}).

 \item  Equation (\ref{state}) is said to be LMPFS with respect  to a subspace $Z$ of  $U$ if  there is
 a $kT$-periodic
$K(\cdot)\in L^\infty\left(\mathbb{R}^+;\cL(H, Z)\right)$ for some  $k\in \mathbb{N}$
 such that the equation (\ref{fstate}) is exponentially stable. Any such a $K(\cdot)$ is called an LMPFS law for Equation (\ref{state}) with respect to  $Z$.
 \end{itemize}

\begin{proposition}\label{lemma3.1}
Let $n_0$ and $H_1$ be given by (\ref{number}) and  (\ref{wgs1.15}) respectively. Suppose that  $Z\subseteq U$ is a finite dimensional subspace satisfying (\ref{3-14}).
 Then, Equation (\ref{state})
is LMPFS with respect to $Z$.
\end{proposition}
\noindent{\it Proof.}
 Let $Z\subseteq U$ satisfy (\ref{3-14}).  We organize the proof by  several steps as follows.

\vskip 5pt
\noindent{\it Step 1. For any $ h_1\in H_1$, to construct control $u^{h_1}(\cdot)\in L^2(\mathbb{R}^+;Z)$ such that $Py(n_0T; 0,h_1,u^{h_1})=0$ }

Because $\mbox{dim} H_1=n_0$ (see (\ref{wgs1.14})), we can set $\{\eta_1,\cdots,\eta_n\}$ to  be an orthonormal basis of $H_1$.
Define a linear map $\mathcal{F}: \mathbb{R}^{n_0}\rightarrow H_1$ by setting
\begin{equation}\label{wgs3.2}
\mathcal{F}(\textbf{a})\deq\sum_{j=1}^{n_0}a_j\eta_j\triangleq(\eta_1,\dots,\eta_{n_0})(a_1,\dots,a_{n_0})^*\;\;\mbox{for each}\;\; \textbf{a}=(a_1,\dots,a_{n_0})^*\in \mathbb{R}^{n_0}.
\end{equation}
Clearly, $\mathcal{F}$ is invertible and
\begin{equation}\label{wgs3.3}
\mathcal{F}^{-1}(h_1)=(\lan h_1,\eta_1 \ran,\dots,\lan h_1,\eta_{n_0}\ran)^*.
\end{equation}
Since $Z$ is finite dimensional, we can assume that $\mbox{dim}Z=m_0<\infty$. Write $\{z_1,\dots, z_{m_0}\}$ for  an orthonormal basis of $Z$. Define a linear map
$\mathcal{G}: L^2(\mathbb{R}^+;\mathbb{R}^{m_0})\rightarrow L^2(\mathbb{R}^+; Z)$ by setting
\begin{equation}\label{wgs3.4}
\mathcal{G}(\beta)(t)\deq(z_1,\dots,z_{m_0})\beta(t)\triangleq\sum_{j=1}^{m_0}\beta_j(t)z_j,\; \mbox{a.e.}\; t\geq 0,
\end{equation}
for each $\beta(\cdot)\triangleq(\beta_1(\cdot),\dots,\beta_{m_0}(\cdot))^*\in L^2(\mathbb{R}^+;\mathbb{R}^{m_0})$ with $\beta_j(\cdot)\in L^2(\mathbb{R}^+;\mathbb{R})$.
Clearly, $\mathcal{G}$ is invertible and it holds that
\begin{equation}\label{wgs3.10}
u(t)=\mathcal{G}\mathcal{G}^{-1}(u)(t)=(z_1,\dots,z_{m_0})\mathcal{G}^{-1}(u)(t)\;\;\mbox{for a.e.}\;\; t\geq 0,\;\;\mbox{when}\;\;u(\cdot)\in L^2(\mathbb{R}^+;Z).
\end{equation}
By the facts that  $ P \cP=\cP P$, $H_1$ is invariant with respect to $\cP$ and $P: H\rightarrow H_1$ is a projection (see Lemma~\ref{lemma2.2}), and by (\ref{2-13}),  we obtain that
$$P \cP^{n_0} h_1=\cP^{n_0} h_1=\cP^{n_0}_1 h_1\;\;\mbox{ for each}\;\; h_1\in H_1.$$
 From this, we see that
 \begin{equation}\label{3-2}\ba{rl}
& Py (n_0T;0,h_1,u) =P \cP^{n_0}h_1+\ds\int^{n_0T}_0P\Phi(n_0T,s)D(s) u(s)ds\\
 &=\cP^{n_0}_1h_1+\ds\int^{n_0T}_0P\Phi(n_0T,s)D(s)u(s)ds\;\;\mbox{for all}\;\;h_1\in H_1, u(\cdot)\in L^2(\mathbb{R}^+;Z).\ea
 \end{equation}
Let $\hat A\in \mathbb{R}^{n_0\times n_0}$ be the matrix of $\cP_1$ under $\{\eta_1,\dots,\eta_{n_0}\}$, i.e.,
$\cP_1(\eta_1,\dots,\eta_{n_0})=(\eta_1,\dots,\eta_{n_0})\hat A$.
Then, it follows from (\ref{wgs3.2}) and (\ref{wgs3.3}) that
\begin{equation}\label{wgs3.8}
\cP_1^{n_0}h_1=\cP_1^{n_0}(\eta_1,\dots,\eta_{n_0})\mathcal{F}^{-1}(h_1)=(\eta_1,\dots,\eta_{n_0})\hat A^{n_0}\mathcal{F}^{-1}(h_1),\;\;\mbox{when}\;\;h_1\in H_1.
\end{equation}
Since $Py (n_0T;0,h_1,u)\in H_1$, it follows by (\ref{wgs3.2}) and (\ref{wgs3.3}) that
\begin{equation}\label{wgs3.9}
Py (n_0T;0,h_1,u)=(\eta_1,\dots,\eta_{n_0})\mathcal{F}^{-1}(Py (n_0T;0,h_1,u))\;\;\mbox{for all}\;\; h_1\in H_1, u(\cdot)\in L^2(\mathbb{R}^+;Z).
\end{equation}
Since $D(\cdot)\in L^\infty(\mathbb{R}^+; \mathcal{L}(U;H))$ and $P\Phi(n_0T,s)D(s)(z_1,\dots,z_{m_0})\in (H_1)^{m_0}$ for a.e. $s\geq 0$, there  is a unique $\hat B(\cdot)\in
L^\infty(\mathbb{R}^+; \mathbb{R}^{n_0\times m_0})$ such that
\begin{equation}\label{wgs3.11}
P\Phi(n_0T,s)D(s)(z_1,\dots,z_{m_0})=(\eta_1,\dots,\eta_{n_0})\hat B(s)\;\;\mbox{for a.e.}\;\; s\in \mathbb{R}^+.
\end{equation}
Now, from (\ref{3-2}), (\ref{wgs3.8}), (\ref{wgs3.9}), (\ref{wgs3.10}) and (\ref{wgs3.11}), we see that for each $h_1\in H_1$ and $u(\cdot)\in L^2(\mathbb{R}^+;Z)$,
\begin{equation}\label{wgs3.12}
\mathcal{F}^{-1}(Py (n_0T;0,h_1,u))=\hat A^{n_0}\mathcal{F}^{-1}(h_1)+\ds\int^{n_0T}_0 \hat B(s)\mathcal{G}^{-1}(u)(s)ds.
\end{equation}
Meanwhile, it follows from (\ref{2-14}), (\ref{csubspace}), (\ref{wgs3.10}), (\ref{wgs3.11}) and  (\ref{wgs3.4}) that
\begin{equation*}
\hat V_{n_0}^Z=\Big\{ (\eta_1,\dots,\eta_{n_0})\ds\int_0^{n_0T}\hat B(s)\beta(s)ds \; \big |\; \beta(\cdot)\in L^2(\mathbb{R}^+; \mathbb{R}^{m_0})\Big\}.
 \end{equation*}
Since $Z$ satisfies (\ref{3-14}), the above equality yields
$$
\mbox{span} \{\eta_1,\dots, \eta_{n_0}\}=\Big\{ (\eta_1,\dots,\eta_{n_0})\ds\int_0^{n_0T}\hat B(s)\beta(s)ds\; \big |\; \beta(\cdot)\in L^2(\mathbb{R}^+; \mathbb{R}^{m_0})\Big\},
$$
which is equivalent to
\begin{equation}\label{wgs3.13}
\Big\{\ds\int_0^{n_0T}\hat B(s)\beta(s)ds\; \big |\; \beta(\cdot)\in L^2(\mathbb{R}^+; \mathbb{R}^{m_0})\Big\}=\mathbb{R}^{n_0}.
\end{equation}
From (\ref{wgs3.13}), we see that the finite-dimensional controlled system $x'(t)=\hat B(t)\beta(t),\; t\geq 0$,
(where $x(\cdot)=(x_1(\cdot),\dots, x_{n_0}(\cdot))^*$  is treated as a state and $\beta(\cdot)$ is treated as a control) is exactly  controllable. Hence, the matrix $\ds\int^{n_0T}_0 \hat B(s)\hat B(s)^*ds$ is positive definite (see, for instance, \cite{Sontag}).

Now for each $h_1\in H_1$, we define  $\beta^{h_1}(\cdot)\in L^2(\mathbb{R}^+; \mathbb{R}^{m_0})$ by setting
\begin{equation}\label{3-4-1}
\beta^{h_1}(t)=\left\{\ba{ll}
\ns-\hat B(t)^*\ds\Big(\int^{n_0T}_0\hat B(s)\hat B(s)^*ds\Big)^{-1}\hat A^{n_0}\mathcal{F}^{-1}(h_1)\;\;&\mbox{for a.e.}\;\; t\in[ 0,n_0T);\\
\ns 0& \mbox{for a.e.}\;\; t\in[n_0T,+\infty).\ea\right.\end{equation}
It is clear that
\begin{equation}\label{3-4-2}
 \hat A^{n_0}\mathcal{F}^{-1}(h_1)+\ds\int_0^{n_0T}\hat B(s)\beta^{h_1}(s)ds=0.
 \end{equation}
Then,  for each $h_1\in H_1$, we construct a control
\begin{equation}\label{3-4}
 u^{h_1}(\cd)=\mathcal{G}(\beta^{h_1})(\cd)\;\;\mbox{over}\;\;\mathbb{R}^+.
\end{equation}
By (\ref{wgs3.4}),  $ u^{h_1}(\cdot)\in L^2(\mathbb{R}^+; Z)$. Meanwhile, it follows from (\ref{3-4-2}) and (\ref{wgs3.12}) that
$$\mathcal{F}^{-1}(Py(n_0T; 0,h_1, u^{h_1}))=0,\;\;\mbox{when}\;\; h_1\in H_1.$$
This, implies that
\begin{equation}\label{wgs3.17}
Py (n_0T;0,h_1, u^{h_1})=0,\;\;\mbox{when}\;\; h_1\in H_1.
\end{equation}

\vskip 5pt

\noindent{\it Step 2. To show the existence of an $N_0\in \mathbb{N}$ such that
\begin{equation}\label{3-6}
\|y(NT;0,h,\cL (P h))\|\le \d_0\|h\|\;\;\;\mbox{for all}\;\; h\in H\;\;\mbox{and}\;\; N\ge N_0,
 \end{equation}
 where $\delta_0\triangleq (1+\bar\delta)/2$ with $\bar\delta$ given by (\ref{decay}). }
\vskip 5pt

Define an  operator $\cL: H_1\rightarrow L^2(\mathbb{R}^+;Z)$ by setting
\begin{equation}\label{wgs3.18}
\cL h_1(\cd)= u^{h_1}(\cd)\;\;\mbox{for all}\;\; h_1\in H_1,
\end{equation}
where $u^{h_1}(\cdot)$ is given by (\ref{3-4}). Several observations on $\cL$ are given in order. First, it is clear that  $\mathcal{L}$ is linear. Second,  from
(\ref{wgs3.18}), (\ref{3-4}), (\ref{3-4-1}) and the fact that $\hat B(\cdot)\in
L^\infty(\mathbb{R}^+; \mathbb{R}^{n_0\times m_0})$, we have
\begin{equation}\label{luis3.19}
\|\cL\|\triangleq\|\cL\|_{\cL(H_1, L^2(\mathbb{R}^+;Z))}<+\infty.
\end{equation}
Next, by (\ref{wgs3.18}), (\ref{3-4}), (\ref{wgs3.4}) and (\ref{3-4-1}), we see that
\begin{equation}\label{luis3.20}
\cL h_1(\cdot)=0\;\;\mbox{over}\;\; [n_0T,+\infty),\;\;\mbox{when}\;\; h_1\in H_1.
\end{equation}
Finally, it follows from  (\ref{wgs3.17}) that
\begin{equation}\label{3-5}
Py(n_0T;0,h_1,\cL h_1)=0,\;\;\mbox{when}\;\;h_1\in H_1.\end{equation}

Let $\rho_0=-\ln\d_0/{T}$. Since $\delta_0\triangleq (1+\bar\delta)/2$, we see  that
$ 0<\rho_0<-\ln\bar \d/{T}\triangleq \bar\rho$.
Then, by Part $(f)$ of  Lemma~\ref{lemma2.2}, there is a constant $C_{\rho_0}>0$ such that
\begin{equation}\label{3-7}
\|y(kT;0,h_2,0)\|=\|\Phi(kT,0)h_2\|\le C_{\rho_0}e^{- \rho_0 kT}\|h_2\|
=C_{\rho_0}\delta_0^k\|h_2\|,\;\;\mbox{when}\;\;k\in\mathbb{N}, h_2\in H_2.
\end{equation}
We  claim that there is a constant $C>0$ such that
\begin{equation}\label{3-9}
\|y(NT;0,h_1,\cL h_1)\|
\le CC_{\rho_0}\delta_0^{N-n_0}\|h_1\|,\;\;\mbox{when}\;\;h_1\in H_1, N\geq n_0,\end{equation}
where $\cL$ is given by (\ref{wgs3.18}).
In fact, because
$$\big \|y(n_0T; 0,h_1,\cL h_1)\big\|\le \ds\big\|\Phi(n_0T,0)h_1\big\|+\left\|\int^{n_0T}_0\Phi(n_0T,s)D(s)\cL
h_1(s)ds\right\|\;\;\mbox{for each}\;\; h_1\in H_1,$$
there is a constant $C>0$ such that
\begin{equation}\label{3-8}
\|y(n_0T;0,h_1,\cL h_1)\|\le C\|h_1\|\;\;\mbox{for all}\;\; h_1\in H_1.
\end{equation}
Here, we used  facts that $D(\cdot)\in L^\infty(\mathbb{R}^+; \mathcal{L}(U;H))$ (see the assumption $(\cH_3)$) and $\cL$ is linear and bounded (see (\ref{luis3.19})).
Meanwhile, it follows from (\ref{luis3.20}), (\ref{wgs1.5}) and (\ref{solution1}) with $u(\cdot)\equiv 0$ that
when $N\geq n_0$ and $h_1\in H_1$,
\begin{equation}\label{wgs3.24}
\ba{rl}\ns&y(NT;0,h_1,\cL h_1)=y(NT;n_0T,y(n_0T;0,h_1,\cL h_1),\cL h_1)\\
\ns=&y(NT;n_0T,y(n_0T;0,h_1,\cL h_1),0)=\Phi(NT, n_0T)y(n_0T; 0,h_1,\cL h_1)\\
\ns=& \Phi((N-n_0)T, 0)y(n_0T; 0,h_1,\cL h_1)=
y((N-n_0)T;0,y(n_0T;0,h_1,\cL h_1),0).\ea
\end{equation}
Because of (\ref{3-5}) and (\ref{2-2})-(\ref{2-2-1}) with $t=0$, it holds that
$y(n_0T;0,h_1,\cL h_1)\in H_2$, when $h_1\in H_1$. This along with
(\ref{wgs3.24}), (\ref{3-7}) and (\ref{3-8}), leads to (\ref{3-9}).

Let
\begin{equation}\label{3-10}
N_0=\max\Big\{  \ds\frac{\ln C_{\rho_0}+\ln\left(C\delta_0^{-n_0}\|P\|+\|I-P\|\right)}{\ln (1/\delta_0)} +2,~n_0\Big\}.
\end{equation}
(Here, $[r]$ with $r\in \mathbb{R}$ denotes the integer such that $r-1< [r]\leq r$.) Then, it follows from (\ref{3-9}), (\ref{3-7}) and (\ref{3-10}) that
\begin{equation*}\ba{rl}
\ns&\|y(NT;0,h, \cL (Ph))\|\leq \|y(NT;0,Ph, \cL (Ph))\|+\|y(NT;0,(I-P)h, 0)\|\\
\ns\leq&
 C_{\rho_0}\delta_0^{N}(C\delta_0^{-n_0}\|P\|+\|I-P\|)\|h\|\le\delta_0 \|h\|,\qq\mbox{when}\;\; N\geq N_0\;\;\mbox{and}\;\; h\in H.
\ea \end{equation*}
This leads to (\ref{3-6}).

\vskip 5pt

\noindent {\it Step 3. To study a  value function associated with  a class of optimal control problems}

\vskip 5pt

Given  $N\in \mathbb{N}$ and  $t\in [0,NT)$, $h\in H$ and $u(\cdot)\in L^2(0,NT;Z)$, consider the equation:
\begin{equation}\label{stateft}
\left\{\ba{ll}y'(s)+A y(s)+B(s)y(s)=D(s)\bigl|_Zu(s)\;\;\mbox{in}\;\; (t,NT), \\
 \ns y(t)=h,\ea\right.
 \end{equation}
where     $D(t)\bigl|_Z$  is the restriction of  $D(t)$  on the subspace $Z$. Because of assumptions $(\cH_1)$-$(\cH_3)$, Equation (\ref{stateft}) has a unique mild solution  $y^Z_N(\cdot; t,h,u)\in C([0,NT];H)$ (see Proposition 5.3 on Page 66 in \cite{Li}).
 Clearly, $y^Z_N(\cdot; t,h,u)=y(\cdot; t,h,u)\big|_{[t,NT]}$.
 For each $\varepsilon>0$, we define a cost functional
$J^{Z}_{\e,N,t,h}(\cdot):L^2(0,NT;Z)\rightarrow \mathbb{R}^+$ by setting
\begin{equation}\label{3-28-1}
J^{Z}_{\e,N,t,h}(u)=\ds\int^{NT}_t\e\|u(s)\|_U^2ds+\|y^Z_N(NT;t,h,u)\|^2,\;\; u\in L^2(0,NT;Z).
\end{equation}
Then, for each  $N\in \mathbb{N}$, $\varepsilon>0$, $t\in [0,NT]$ and  $h\in H$, we define the optimal control problem
\begin{equation*}\label{Luis3.30}
(P)^{Z}_{\e, N,t,h}:\;\;\;\;\;\;\; \ds\inf_{u\in L^2(0,NT;\,Z) }J^{Z}_{\e,  N, t,h}(u).
\end{equation*}
\noindent This a  classical
linear quadratic optimal control problem (see  Page 370 in \cite{Li}). For each $\varepsilon>0$ and $N\in \mathbb{N}$, the value function associated with the above optimal control problems
is
\begin{equation}\label{WGS3.28}
W^{Z}_{\e,N}(t,h)\deq\ds\inf_{u\in L^2(0,NT;\,Z) }J^{Z}_{\e, N, t,h}(u),\;\; t\in [0,NT]\;\;\mbox{and}\;\; h\in H.
\end{equation}
Let
\begin{equation}\label{wgs3.28}
\e_0\deq (\delta_0-\delta_0^2)/{(\|\cL\|\|P\|+1)^2},
 \end{equation}
 where $\delta_0\in (0,1)$ and $\cL$ are given by (\ref{3-6}) and (\ref{wgs3.18}) respectively.
Because of (\ref{luis3.19}), it holds that $0<\e_0<+\infty$.
We claim  that
\begin{equation}\label{wgs3.29}
W^{Z}_{\e, N}(0, h)\le \delta_0\|h\|^2\;\mbox{ for all }\;h\in H,\;\;\mbox{when}\;\; N\ge N_0\;\;\mbox{and}\;\;\e\in (0, \e_0],
\end{equation}
where $N_0$ is given by (\ref{3-10}).
In fact, it follows from (\ref{wgs3.28}) that
\begin{equation}\label{luis3.31}
\e\|\cL (Ph)(\cdot)\|_{L^2(\mathbb{R}^+;Z)}^2\le \e_0\| \cL\|^2\|P\|^2 \|h\|^2 \le (\delta_0-\delta_0^2)\|h\|^2,\;\;\mbox{when}\;\; h\in H\;\;\mbox{and}\;\;\e\in (0,\e_0].
\end{equation}
By (\ref{wgs3.18}) and the fact that $P$ is a projection from $H$ to $H_1$ (see Lemma~\ref{lemma2.2}), we find
\begin{equation}\label{luis3.32}
\cL(Ph)\in L^2(\mathbb{R}^+;Z)\;\;\mbox{for all}\;\; h\in H.
\end{equation}
Since $y^Z_N(\cdot; t,h,u\bigl|_{(0,NT)})=y(\cdot; t,h,u)\big|_{[t,NT]}$ for any $u\in L^2(\mathbb{R}^+;Z)$, it holds that
\begin{equation*}\label{luis3.36}
y^Z_N(NT; 0,h,\cL(Ph)\bigl|_{(0,NT)})=y(NT; 0,h, \cL(Ph)).
\end{equation*}
This, together with  (\ref{WGS3.28}), (\ref{3-28-1}), (\ref{luis3.32}), (\ref{luis3.31}) and (\ref{3-6}), indicates that
\begin{equation*}\label{wgs3.31}
W^{Z}_{\e, N}(0, h)\leq J^{Z}_{\e, N, 0,h}(\cL(Ph)\bigl|_{(0,NT)})=\e\ds\int_0^{NT}\|\cL(Ph)(s)\|^2ds+\|y^Z_N(NT;0,h,\cL(Ph))\|^2\le \delta_0\|h\|^2,
\end{equation*}
 when  $N\ge N_0$, $\e\in (0,\e_0]$, $h\in H$, i.e.,  (\ref{wgs3.29}) stands.

\vskip 5pt

\noindent {\it Step 4. To construct an  $NT$-periodic $K^Z_{\varepsilon,N}(\cdot)\in L^\infty(\mathbb{R}^+; \cL(H,Z))$ }

Arbitrarily fix an  $\e\in (0,\e_0]$ and an $N\geq N_0$, where $N_0$ and $\e_0$ are given by (\ref{3-10}) and (\ref{wgs3.28}) respectively.
 By the exactly same way to show Corollary 2.10 on Page 379 and   Theorem 4.3 on Page 397 in \cite{Li}, we can verify that
$$W^{Z}_{\e,N}(t,h)=\lan Q^{Z}_{\e,N}(t)h,h\ran,\;\;\mbox{when}\;\; h\in H,$$
where  $Q^{Z}_{\e,N}(\cdot)\in C([0,NT);\cL(H))$ has the following properties:
 $(i)$ for each $t\geq 0$, $Q^{Z}_{\e,N}(t)$ is self-adjoint;   $(ii)$  it solves  the
 Riccati integral equation:
\begin{equation}\label{wgx3.37}
 \ba{r}
 {Q}(t)h=\Phi(NT,t)^*\Phi(NT,t)h
 -\ds\frac{1}{\e}\int^{NT}_t \Phi(s,t)^* Q(s)^* D(s)\bigl|_Z\left(D(s)\bigl|_Z\right)^*Q(s)\Phi(s,t)hds,\\
 {\mbox{ for all }} h\in H{\mbox{ and }}t\in[0,NT).\ea
\end{equation}
Besides, it follows form  (\ref{3-28-1})  that
$$0\le\left\langle h, Q^{Z}_{\e,N}(t)h\right\rangle\le J^Z_{\e, N,t,h}(0)\le\|\Phi(NT,t)\|^2\|h\|^2 \;\mbox{ for all }\;h\in H.$$
Define $K^{Z}_{\e,N}(\cdot): [0,NT)\rightarrow \cL(H;Z)$ by
\begin{equation}\label{3-12}
 K^{Z}_{\e,N}(t)=-\ds\frac{1}{\e}\left(D(s)\bigl|_Z\right)^* Q^{Z}_{\e,N}(t)\;\;\mbox{for a.e.}\;\; t\in [0,NT).
\end{equation}
One can easily check that
$
K^{Z}_{\e,N}(\cd)\in  L^\infty(0,NT;\cL(H;Z))$.
From this and  assumptions $(\cH_1)$-$(\cH_3)$, we see that (see Proposition 5.3 on Page 66 in \cite{Li})
 the feedback equation:
\begin{equation}\label{3-12-4} \left\{\begin{array}{lll}
\ns y'(s)+Ay(s)+B(s)y(s)= D(s)\bigl|_Z K^{ Z}_{\e,N}(s)y(s)\;\;&\mbox{in}\;\; (0,NT),\\
\ns y(0)=h\in H
\end{array}\right.\end{equation}
has a unique mild solution   $y^Z_{\e,N}(\cdot; 0,h)\in C([0,NT];H)$.
Let
\begin{equation}\label{3-12-5}
u^Z_{\e,N,0,h}(s)\deq K^Z_{\e,N}(s)y^Z_{\e,N}(s; 0, h)\qq \mbox{for a.e.}\;\; s\in(0, NT).\end{equation}
By the state feedback representation of optimal controls for linear quadratic control problems (see Section 3.4 in Chapter 9 in \cite{Li}, in particular,  (3.71) on Page 392 and (4.12) on Page 397 in \cite{Li}), we see that
$u^{Z}_{\e,N,0,h}(\cd)$  defined by (\ref{3-12-5}) is the optimal control to $(P)^Z_{\varepsilon,N,0,h}$.
This, along with  (\ref{WGS3.28}), yields that
\begin{equation}\label{luis3.42}
W^Z_{\e, N}(0,h)=J^Z_{\e, N, 0,h}\big( u^{Z}_{\e,N,0,h}\big),\;\;\mbox{when}\;\; h\in H.
\end{equation}
By (\ref{stateft}) with $t=0$, (\ref{3-12-4}) and (\ref{3-12-5}), we see that
$y^Z_{\e,N}(NT; 0,h)=y^Z_N(NT; 0,h, u^Z_{\e,N,0,h})$. From  this, (\ref{3-28-1}), (\ref{luis3.42}) and  (\ref{wgs3.29}), it follows that
\begin{equation}\label{3-12-6}
\|y^{Z}_{\e,N}(NT; 0,h)\|^2\le J^Z_{\e,N,0,h}(u^Z_{\e,N,0,h})=W^Z_{\e, N}(0,h)\le\d_0\|h\|^2,\;\;\mbox{when}\;\; h\in H.
\end{equation}
Now, we extend $NT$-periodically $K^{Z}_{\varepsilon,N}(\cdot)$ over $\mathbb{R}^+$ by setting
\begin{equation}\label{3-12-1}
 K^{Z}_{\varepsilon,N}(t+kNT)= K^{Z}_{\varepsilon,N}(t)\;\;\mbox{for all}\;\; t\in [0,NT),~k\in\mathbb{N}.
\end{equation}
Clearly, $ K^{Z}_{\e,N}(\cd)\in  L^\infty(\mathbb{R}^+;\cL(H;Z))$ is $NT$-periodic.

\vskip 5pt

\noindent {\it Step 5. To prove that when $\e\in (0,\e_0]$ and $N\geq N_0$ (where $N_0$ and $\e_0$ are given by (\ref{3-10}) and (\ref{wgs3.28}), respectively), $K^Z_{\e,N}(\cdot)$ defined by (\ref{3-12}) and (\ref{3-12-1}) is an LMPFS law for Equation (\ref{state}) with respect to $Z$}

\vskip 5pt

Consider the feedback equation:
\begin{equation}\label{3-13} \left\{\begin{array}{lll}
\ns y'(s)+Ay(s)+B(s)y(s)= D(s)\bigl|_Z K^{ Z}_{\e,N}(s)y(s)\qq&\mbox{in}\;\; \mathbb{R}^+,\\
\ns y(0)=h\in H.
\end{array}\right.\end{equation}
By  assumptions $(\cH_2)$ and $(\cH_3)$, and by the fact that $ K^{Z}_{\e,N}(\cd)\in  L^\infty(\mathbb{R}^+;\cL(H;Z))$, we have
$$B(\cdot)-D(\cdot)\bigl|_Z K^{Z}_{\e,N}(\cdot)\in L^1_{loc}(\mathbb{R}^+;\cL(H)).$$
Thus,  for each $h\in H$, Equation (\ref{3-13}) has a unique mild solution  $y^Z_\e(\cdot;0,h)\in C(\mathbb{R}^+; H)$ (see  Proposition 5.3 on Page 66 in \cite{Li}).
 Clearly,
$$
y^Z_\e(t;0,h)=y^Z_{\e, N}(t;0,h)\;\;\mbox{for each}\;\; t\in[0,NT].
$$
 Write
 $\{\Phi^Z_{\e, N}(t, s)\bigl\}_{0\le s\le t<+\infty}$ for  the evolution system generated by
$-A-B(\cd)+D(\cdot)\bigl|_Z K^Z_{\e,N}(\cdot)$.
Then
it holds that (see Proposition 5.7, Page 69, \cite{Li})
\begin{equation*}\label{wgs3.36}
y^Z_{\e,N}(NT; 0,h)=y^Z_\e(NT;0,h)=\Phi^Z_{\e,N}(NT,0)h\;\;\mbox{ for all }\;\;h\in H.
\end{equation*}
This, along with (\ref{3-12-6}) and the fact that $\delta_0<1$ (see (\ref{3-6})), yields
\begin{equation}\label{wgs3.37}
\|\Phi^Z_{\e,N}(NT,0)\|\le \sqrt{\delta_0} <1.
\end{equation}
Since $B(\cdot)$ and $D(\cdot)$ are $T$-periodic and
 $K^Z_{\e,N}(\cdot)$ is $NT$-periodic, it follows that
\begin{equation}\label{WHUANGLIN3.38}
\Phi^Z_{\e,N}(t+NT,s+NT)= \Phi^Z_{\e,N}(t,s),\;\;\mbox{when}\;\;0\leq s\leq t<+\infty.
\end{equation}
By (\ref{WHUANGLIN3.38}) and  (\ref{wgs3.37}), one can easily  shows that  Equation (\ref{3-13}) is exponentially
 stable.
 Hence,   $K^Z_{\e,N}(\cdot)$, with $\e\in (0,\e_0]$ and $N\geq N_0$, is an LMPFS law for Equation (\ref{3-13}). This completes the proof. \endpf

\begin{proposition}\label{lemma3.3} Let   $Z$ be   a subspace of $U$. Then, Equation (\ref{state})  is LPFS with respect to $Z$ if and only if
it is LMPFS with respect to $Z$.
\end{proposition}
{\it Proof.} It is clear that
$\mbox{Equation}\;(\ref{state})\;\mbox{is LPFS w.r.t.}\;Z\Rightarrow \mbox{Equation}\;(\ref{state})\;\mbox{is LMPFS w.r.t.}\;Z$.
We will show the reverse.
Suppose that  Equation (\ref{state}) is LMPFS with respect to $Z$. Then there are an $N\in \mathbb{N}$ and an $NT$-periodic $\hat K^Z_N(\cdot)\in L^\infty(\mathbb{R}^+;\mathcal{L}(H;Z))$
such that the feedback equation
\begin{equation}\label{3-100}
y'(s)+Ay(s)+B(s)y(s)= D(s)\bigl|_Z \hat K^Z_N(s)y(s),\;\; s\geq 0
\end{equation}
is exponentially stable. From this, assumptions $(\mathcal{H}_1)$-$(\mathcal{H}_3)$ and the fact that $\hat K^Z_N(\cdot)\in L^\infty(\mathbb{R}^+;\mathcal{L}(H;Z))$ is $NT$-periodic, one can easily verify that
for each $t\geq 0$ and $h\in H$,
the solution $\hat y^Z_N(\cdot; t, h)$ to the equation
\begin{equation}\label{3-20-1}\left\{\ba{l}
y'(s)+Ay(s)+B(s)y(s)= D(s)\bigl|_Z \hat K^Z_N(s)y(s) \;\;\mbox{in}\;\; [t,\infty),\\
y(t)=h\ea\right.
\end{equation}
satisfies
\begin{equation}\label{wgs3.49}
\|\hat y^Z_N(s;t,h)\|\leq M_1e^{-\d_1 (s-t)}\|h\|,\;\;\mbox{when}\;\; s\geq t\;\;\mbox{and}\;\; h\in H,
\end{equation}
where $M_1$ and $\delta_1$ are two positive constants independent of $h$, $t$ and $s$. Write
\begin{equation}\label{wgs3.50}
\hat C_0\triangleq \left\| \hat K^Z_N (\cd)\right\|_{L^\infty(\mathbb{R}^+;\cL(H;Z))}.
\end{equation}
The rest of the  proof is organized by two steps as follows.

\vskip 5pt
\noindent {\it Step 1. To study a value function associated with a class of optimal control problems}

\vskip 5pt

For each $t\geq 0$ and  $h\in H$, we define an optimal control problem:
\begin{equation}\label{wgs3.51}
 (P)^Z_{t,h}:\;\;\;\;\;\ds\inf\limits_{u(\cdot)\in L^2(\mathbb{R}^+;Z)}\int^{\infty}_t\left(\|y(s;t,h,u)\|^2+ \|u(s)\|_U^2\right)ds.
\end{equation}
Associated with $\big\{(P)^Z_{t,h}\big\}_{t\geq 0, h\in H}$, we define  a value function $W^Z(\cdot,\cdot):\mathbb{R}^+\times H\rightarrow \mathbb{R}$ by
\begin{equation}\label{3.20}
W^Z(t, h)=\ds\inf\limits_{u(\cdot)\in L^2(\mathbb{R}^+;Z)}\int^{\infty}_t\left(\|y(s;t,h,u)\|^2+ \|u(s)\|_U^2\right)ds, \;\; (t,h)\in\mathbb{R}^+\times H.
\end{equation}
It is well defined. In fact, define a control $u^Z_{N, t, h}(\cdot)$ by setting
\begin{equation*}
u^Z_{N, t, h}(s)\triangleq\left\{\ba{ll}
 K^Z_N(s)\hat y^Z_{N}(s; t, h)\;\;&\mbox{for a.e.}\;\; s\in [t, \infty),\\
\ns 0& \mbox{for a.e.}\;\; s\in[0, t).\ea\right.
\end{equation*}
By (\ref{wgs3.49}), we see that $u^Z_{N, t, h}(\cdot)\in L^2(\mathbb{R}^+;Z)$.
Meanwhile, it is clear that $\hat y^Z_N(\cdot; t,h)=y(\cdot; t,h, u^Z_{N, t, h})$. These, along with (\ref{wgs3.49}) and (\ref{wgs3.50}), yields that for each $t\geq 0$ and $h\in H$,
\begin{equation}\label{WUHENHEN5.53}
\ba{rl}
\ns&0\leq W^Z(t, h)\le\ds
\int^{\infty}_t\left(\big\|\hat y^Z_N(s;t,h)\big\|^2+ \big\| u^Z_{N, t, h}(s))\big\|_U^2\right)ds\\
\ns&\le\ds
\int^{\infty}_t\big(1+ \big\| \hat K^Z_N (s)\big\|_{\cL(H;Z)}^2\big)\big\|\hat y^Z_N(s;t,h)\big\|^2ds\le\ds\frac{(1+\hat C_0^2)M^2_1}{2\d_1}\big\|h\big\|^2.
\ea
\end{equation}
Thus $W^Z(\cdot,\cdot)$ is well-defined.
 By (\ref{3.20}), one can directly check that when $h,g\in H, t\geq 0$, $ c\in \mathbb{R}$,
$$W^Z(t, ch)=c^2W^Z(t,h)\;\;\mbox{and}\;\; \; W^Z(t,h+g)+W^Z(t, h-g)=2W^Z(t,h)+2W^Z(t,g).
$$
 These, together with (\ref {WUHENHEN5.53}), imply that (see \cite{Kurepa}) there is a unique $Q^Z(\cdot): \mathbb{R}^+\rightarrow\mathcal{L}(H)$,
 with $Q^Z(t)$ self-adjoint for each $t\geq 0$, such that
\begin{equation}\label{WUHENHEN5.54}
W^Z(t,h)=\lan Q^Z(t)h,h\ran\;\;\mbox{for all}\;\; t\in \mathbb{R}^+\;\;\mbox{and}\;\;h\in H.
\end{equation}
This, together with (\ref{WUHENHEN5.53}), implies that
$$0\le Q^Z(t)\le \ds\frac{(1+\hat C_0^2)M^2_1}{2\d_1}I\;\;\mbox{for all}\;\; t\in \mathbb{R}^+.$$
Meanwhile, from the $T$-periodicity of $B(\cdot)$ and $D(\cdot)$, one can easily derive the $T$-periodicity of $W^Z(\cdot,h)$ for each $h\in H$. Thus, by (\ref{WUHENHEN5.54}), $Q^Z(\cdot)$ is $T$-periodic, i.e.,
\begin{equation}\label{WHUANGHEN5.55}
Q^Z(t)=Q^Z(T+t)\;\;\mbox{for all}\;\; t\in \mathbb{R}^+.
\end{equation}
Now we  present other properties of $Q^Z(\cd)$.
By the Bellman optimality principle (see Section 1, Chapter 6 in \cite{Li}), it holds that for any  $t\in [0,T]$ and $h\in H$,
\begin{equation}\label{Valued}
W^Z(t,h)=\ds\inf\limits_{u(\cdot)\in L^2(\mathbb{R}^+;Z)}\Big\{\int^{T}_t\left(\|y(s;t,h,u)\|^2+ \|u(s)\|_U^2\right)ds+W^Z(T,y(T;t,h,u))\Big\}.
\end{equation}
 This, along with (\ref{WUHENHEN5.54}) and (\ref{WHUANGHEN5.55}), yields that for any  $t\in [0,T]$ and $h\in H$,
\begin{equation}\label{3-53-1}
W^Z(t,h)=\ds\inf\limits_{u\in L^2(\mathbb{R}^+;Z)}\Big\{\int^{T}_t\left(\|y(s;t,h,u)\|^2+ \|u(s)\|_U^2\right)ds+\|\big({Q^Z(0)}\big)^{1/2} y(T;t,h,u)\|\Big\},
\end{equation}
i.e.,  $W^Z(\cdot,\cdot)\big|_{[0,T]\times H}$ is the value function associated with the  LQ problems $(P)^{Z,T}_{t,h}$, $t\in [0,T]$, $h\in H$:
\begin{equation}\label{wgs3.57}
\ds\inf\limits_{v\in L^2(0,T;Z)}\Big\{\int^{T}_t\left(\|y^T(s;t,h,v)\|^2+ \|v(s)\|_U^2\right)ds+\lan y^T(T;t,h,v), Q^Z(0)y^T(T;t,h,v)\ran\Big\},
\end{equation}
where $y^T(\cdot; t,h,v)$ is the solution of Equation (\ref{state}) (over $[t,T]$), with the initial condition that $y(t)=h$ and with the  control $v(\cdot)\in L^2(0,T; Z)$.
Furthermore,
by the exactly same way to show Corollary 2.10 on Page 379 and   Theorem 4.3 on Page 397 in \cite{Li},  we can obtain that
\begin{equation}\label{wgs3.58}
W^Z(t,h)=\lan h, \hat Q(t)h\ran,\;\;\mbox{when}\;\; t\in [0,T]\;\;\mbox{and}\;\; h\in H,
\end{equation}
where $\hat Q(\cdot)\in C([0,T);\cL(H))$, with $\hat Q(t)$ self-adjoint and non-negative for each $t\in [0,T]$, satisfies the
following Riccati integral equation:
\begin{equation}\label{wgs3.59}
\ba{rl}
\hat Q(t)h=&\Phi(T,t)^*Q(0)\Phi(T,t)h+\ds\int^T_t\Phi(s,t)^*\Phi(s,t)hds \\
 \ns&-\ds\int^T_t \Phi(s,t)^*\hat Q(s)^* D(s)\bigl|_Z\left(D(s)\bigl|_Z\right)^*\hat Q(s)\Phi(s,t)hds\;\;\mbox{for all}\;\;  h\in H, t\in[0,T).
  \ea
\end{equation}
Because  $Q(t)$ and $\hat Q(t)$ are selfadjoint,
 it follows from (\ref{WUHENHEN5.54}) and (\ref{wgs3.58}) that
  \begin{equation}\label{3-57-1}
  Q^Z(\cdot)\bigr|_{[0,T]}=\hat Q(\cdot).
  \end{equation}

\medskip

\noindent{\it Step 2. To construct a $T$-periodic $\widetilde{K}(\cdot)\in L^\infty(\mathbb{R}^+;\cL(H;Z))$}

Let $Q^Z(\cdot)$ be given by  (\ref{WUHENHEN5.54}). We define  $\widetilde{K}(\cdot)\in L^\infty(\mathbb{R}^+;\cL(H;Z))$ by setting
\begin{equation}\label{3-58-1}
\widetilde{K}(s)=-\ds\left(D(s)\bigl|_Z\right)^*Q^Z(s)\;\;\mbox{for a.e.}\;\; s\in\mathbb{R}^+.
\end{equation}
Since both $Q^Z(\cdot)$ and $D(\cdot)$ are $T$-periodic, so is  $\widetilde{K}(\cdot)$.
Let $\widetilde{y}_T(\cd;0,g)\in C([0,T]; H)$ be the  mild solution to the feedback equation:
\begin{equation*}
\left\{\ba{l}
\ns y'(s)+Ay(s)+B(s)y(s)=D(s)\bigl|_Z\widetilde{K} (s) y(s)\;\;\mbox{for a.e.}\;\; s\in(0,T);\\
\ns y(0)=g\in H.\ea\right.
\end{equation*}
The existence and uniqueness of the solution to the above equation is ensured by of assumptions $(\cH_1)$-$(\cH_3)$ and the fact that $\widetilde{K}(\cdot)\in L^\infty(\mathbb{R}^+;\cL(H;Z))$
(see Proposition 5.3 on Page 66 in \cite{Li}).
For each $g\in H$, we define
\begin{equation*}\label{3-60-1}
 \widetilde{u}_g(t)\triangleq \widetilde{K}(t)\widetilde{y}_T(t;0,g)\;\;\mbox{for a.e.}\;\; t\in(0,T).\end{equation*}
By the state feedback representation of optimal controls for linear quadratic control problems (see Section 3.4 in Chapter 9, in particular the formula (3.71) on Page 392 in \cite{Li}),
 it follows from (\ref{3-58-1}) and (\ref{3-57-1}) that
$ \widetilde{u}_{g}(\cd)$ is the optimal control to $(P)^{Z,T}_{0, g}$ (which is defined by (\ref{wgs3.57})).
By (\ref{3-53-1}), (\ref{wgs3.58}) and (\ref{3-57-1}),  we see that for each $g\in H$,
\begin{equation}\label{3-59-1}
\lan g, Q^Z(0)g\ran=\ds\int^T_0 \|\widetilde{y}_T(t; 0,g)\|^2+\|\widetilde{u}_{g}(t)\|_U^2dt+\lan \widetilde{y}_T(T;0,g), Q^Z(0) \widetilde{y}_T(T;0,g)\ran.
\end{equation}
 Associated with each $h\in H$, we define three sequences $\{h_{k}\}_{k=0}^\infty\subseteq H$, $\{y^h_{k}(\cd)\}_{k=1}^\infty\subseteq C([0,T];H)$ and
 $\{u^h_{k}(\cd)\}_{k=1}^\infty\subseteq L^2(0,T;Z)$ by setting
\begin{equation}\label{3-59-3}
 h_{0}\triangleq h,\; h_{k}\triangleq\widetilde{y}_T(T; 0,h_{k-1}), \; y^h_{k}(\cd)\triangleq\widetilde{y}_T(\cd; 0,h_{k-1}),\; u^h_{k}(\cd)\triangleq\widetilde{u}_{h_{k-1}}(\cd)\;\;\mbox{for all}\;\; k\in\mathbb{N},\end{equation}
Taking $g=h_{k-1}$ in (\ref{3-59-1}), we find that
\begin{equation}\label{3-59-2}
\lan h_{k-1}, Q^Z(0)h_{k-1}\ran=\ds\int^T_0 \|y^h_{k}(t)\|^2+\|u^h_{k}(t)\|_U^2dt+\lan h_{k} , Q^Z(0) h_{k}\ran\;\;\mbox{for each}\;\; k\in\mathbb{N}.
\end{equation}

\vskip 5pt

\noindent {\it Step 3. To prove that $\widetilde{K}$ given by  (\ref{3-58-1})
 is  an LPFS law for Equation (\ref{state}) with respect to  $Z$}

\vskip 5pt

\noindent
Consider the feedback   equation:
\begin{equation}\label{Pfeedback}
y'(t)+Ay(t)+B(t)y(t)=D(t)\bigl|_Z\widetilde{K} (t) y(t)\;\;\mbox{in}\;\; \mathbb{R}^+.
\end{equation}
Because of assumptions $(\cH_1)$-$(\cH_3)$ and the fact that $\widetilde{K}(\cdot)\in L^\infty(\mathbb{R}^+;\cL(H;Z))$, corresponding to each initial condition that $y(0)=h\in H$, Equation (\ref{Pfeedback})  has a unique mild solution  $y_{\widetilde{K}}(\cdot; 0,h)\in C(\mathbb{R}^+; H)$
(see Proposition 5.3 on Page 66 in \cite{Li}).
Let $\{\Phi_{\widetilde{K}}(t,s)\}_{t\geq s\geq 0}$
for the evolution system generated by $-A-B(\cdot)+D(\cdot)\bigl|_Z\widetilde{K} (\cdot)$. Then  $y_{\widetilde{K}}(t; 0,h)=\Phi_{\widetilde{K}}(t,0)h$ for each $t\geq 0$ and  $h\in H$ (see  Proposition 5.7  on Page 69 in \cite{Li}).
By the $T$-periodicity of $B(\cdot)$, $D(\cdot)$ and $\widetilde{K}(\cdot)$, and by (\ref{3-59-3}), one can easily check that
\begin{equation}\label{wgs3.64}
\Phi_{\widetilde{K}}(t+T,s+T)=\Phi_{\widetilde{K}}(t,s),\;\;\mbox{when}\;\; t\ge s\geq 0.
\end{equation}
Meanwhile, it follows from the definition of $y^h_k$ (see (\ref{3-59-3})) that given $h\in H$,
\begin{equation}\label{3-66-1}
y_{\widetilde{K}}(t; 0,h)=y^h_{[t/T]+1}(t-[t/T]T)\;\;{\mbox{ for all }}\;\;t\in\mathbb{R}^+.
\end{equation}
 We first claim that
\begin{equation}\label{3-23}
\int^\infty_0\| \Phi_{\widetilde{K}}(t,0)h\|^2dt\le \|Q^Z(0)\|\|h\|^2,\;\;\mbox{when}\;\; h\in  H.
\end{equation}
Indeed, it follows from (\ref{3-66-1}) that
\begin{equation}\label{3-68-1}
\int^\infty_0\| \Phi_{\widetilde{K}}(t,0)h\|^2dt
=\sum\limits^\infty_{k=1}\int^{kT}_{(k-1)T}\|  y_{\widetilde{K}}(t,0,h)\|^2dt=\sum\limits^\infty_{k=1}\int^{T}_0\| y^h_{k}(t)\|^2dt.\end{equation}
Meanwhile, by (\ref{3-59-2}), we find that
$$
\sum\limits^N_{k=1}\lan h_{k-1}, Q^Z(0)h_{k-1}\ran=\sum\limits^N_{k=1}\big[\ds\int^T_0 \|y^h_{k}(t)\|^2+\|u^h_{k}(t)\|_U^2dt+\lan h_{k} , Q^Z(0) h_{k}\ran\big],\qq N\in\mathbb{N}.
$$
Thus
$$\ba{rl}\ns&\ds\sum\limits^\infty_{k=1}\int^{T}_0\| y^h_{k}(t)\|^2dt
=\lim\limits_{N\rightarrow\infty}\sum\limits^N_{k=1}\ds\int^T_0 \|y^h_{k}(t)\|^2dt\\
=&\lim\limits_{N\rightarrow\infty}\big[\ds\lan h, Q^Z(0)h\ran-\sum\limits^N_{k=1} \ds\int^T_0  \|u^h_{k}(t)\|_U^2dt-\lan h_{N} , Q^Z(0) h_{N}\ran\big]\\
\le&\lan h, Q^Z(0)h\ran\le\|Q^Z(0)\|\|h\|^2.
\ea
$$
This, along with (\ref{3-68-1}), leads to (\ref{3-23}).

Since $\widetilde{K}(\cdot)\in L^\infty(\mathbb{R}^+;\cL(H;U))$ is T-periodic, in order to show that
$\widetilde{K} (\cdot)$ is an LPFS law for Equation (\ref{state}) with respect to $Z$, it suffices to prove that there are positive $M$ and $\delta$ such that $\|y_{\tilde K(t; 0,h)}\|\leq Me^{-\delta t}\|h\|$ for all $h\in H$.
This will be done if one can show that $\hat\delta<1$,
where
\begin{equation}\label{wgs3.69}
\hat\delta\triangleq \lim\limits_{k\rightarrow\infty}\|\Phi_{\widetilde{K}}(T,0)^k\|^{1/k}.
\end{equation}
The reason is that  $\Phi_{\tilde K}(\cdot,\cdot)$ is $T$-periodic and  $y_{\widetilde{K}}(t; 0,h)=\Phi_{\widetilde{K}}(t,0)h$.

The rest is to prove  that $\hat\delta<1$. First we can assume that
\begin{equation}\label{wgs3.71}
\hat\rho\triangleq{\ln(\hat\delta+1)}/{T}>0,
\end{equation}
for otherwise $\hat\delta=0<1$. By (\ref{wgs3.69}) and (\ref{wgs3.71}), one can easily check that there is a constant $\hat C_1>0$ such that
\begin{equation}\label{3-23-1}
\left\|\Phi_{\widetilde{K}}(kT,0)\right\|\le \hat C_1e^{\hat\rho kT}\;\;\mbox{for all}\;\; k\in\mathbb{N}.
\end{equation}
Because of $(\mathcal{H}_1)$, $(\mathcal{H}_2)$, $(\mathcal{H}_3)$ and the fact that $\widetilde{K}(\cdot)\in L^\infty(\mathbb{R}^+;\cL(H;U))$, one can easily check that
 $\{\|\Phi_{\widetilde{K}}(t,s)\|_{\mathcal{L}(H)}\}_{0\leq s\leq t\leq T}$ is bounded. Thus,  we can write
\begin{equation}\label{wgs3.73}
\hat C_2\triangleq\max\limits_{0\le t_1\le t_2\le T}\left\|\Phi_{\widetilde{K}}(t_2,t_1)\right\|\in \mathbb{R}^+\;\;\mbox{and}\;\; \hat C_3\triangleq\max\{\hat C_2, ~\hat C_1\hat C_2^2\}\in \mathbb{R}^+.
\end{equation}
We claim that
\begin{equation}\label{3-24}
\| \Phi_{\widetilde{K}}(t,s)\|\le \hat C_3e^{\hat\rho (t-s)}, \;\;\mbox{when}\;\; t\ge s\ge0.
\end{equation}
In fact, given $0\le s\le t$, there are only two cases: $(i)$
 $\left[{t}/{T}\right]=\left[{s}/{T}\right]$ and $(ii)$ $\left[{t}/{T}\right]\neq\left[{s}/{T}\right]$.
 In the first case, we have
 $$ 0\leq t-\left[{s}/{T}\right]\leq T\;\;\mbox{and}\;\; 0\leq s-\left[{s}/{T}\right]\leq T.$$
 These, along with (\ref{wgs3.64}) and (\ref{wgs3.73}), yields
 \begin{equation*}
 \|\Phi_{\widetilde{K}}(t,s)\|=\|\Phi_{\widetilde{K}}(t-\left[{s}/{T}\right]T,s-\left[{s}/{T}\right]T)\leq\hat C_2\leq \hat C_3e^{\hat\rho(t-s)},
 \end{equation*}
 i.e., (\ref{3-24}) holds for the first case. In the second case, we have
 $$
 \left[t/T\right]T\geq\left [s/T\right]T+T\;\;\mbox{and}\;\;(\left[t/T\right]-\left[s/T\right]-1)T\leq t-s.
 $$
These, together with  (\ref{wgs3.64}), (\ref{3-23-1}) and (\ref{wgs3.73}), yields
$$\ba{rl}
&\big\|\Phi_{\widetilde{K}}(t,s)\big\|\leq\big\|\Phi_{\widetilde{K}}\big(t,\,\big[{t}/{T}\big]T\big)\big\|\cdot
\big\|\Phi_{\widetilde{K}}\big(\big[{t}/{T}\big]T,\,\big[{s}/{T}\big]T+T\big)\big\|\cdot
\big\|\Phi_{\widetilde{K}}\big(\big[{s}/{T}\big]T+T,\,s\big)\big\|\\
\ns&\le\big\|\Phi_{\widetilde{K}}\big(t-\big[{t}/{T}\big]T,\,0\big)\big\|\cdot
\big\|\Phi_{\widetilde{K}}\big(\big(\big[{t}/{T}\big]-\big[{s}/{T}\big]-1\big)T,~0\big)\big\|\cdot
\big\|\Phi_{\widetilde{K}}\big(T,\,s-\big[{s}/{T}\big]T\big)\big\|\\
\ns&\le\hat C_1\hat C_2^2e^{\hat\rho (t-s)}\le \hat C_3 e^{\hat\rho (t-s)},
\ea$$
%
i.e., (\ref{3-24}) holds for the second case. In summary, we conclude that (\ref{3-24}) stands.

Let
\begin{equation}\label{WANGwgs3.76}
\hat C_4\triangleq\max\Big \{\max_{t\in[0,T]}\|\Phi_{\widetilde{K}}(t,0)\|, \; \hat C_3 \sqrt{\|Q^Z(0)\|{2\hat\rho}/(1-e^{-2\hat\rho T}}) \Big\}\;\;\mbox{and}\;\; \hat C_5\triangleq\max\{\hat C_2, \hat C_4\hat C_2\},
\end{equation}
where $\hat C_2$ and $\hat C_3$ are given by (\ref{wgs3.73}), $\hat\rho$ is given by (\ref{wgs3.71}).
Then, we claim
\begin{equation}\label{3-24-2}
\| \Phi_{\widetilde{K}}(t,s)\|\le \hat C_5,\;\;\mbox{when}\;\; t\geq s\geq 0 \;(\mbox{where}\; \hat C_5\; \mbox{is given by}\; (\ref{WANGwgs3.76})).
\end{equation}
For this purpose, we first show
\begin{equation}\label{wgs3.77-1}
\|\Phi_{\widetilde{K}}(t,0)\|\leq \hat C_4\;\;\mbox{for all}\;\; t\geq 0 \;(\mbox{where}\; \hat C_4\; \mbox{is given by}\; (\ref{WANGwgs3.76})).
\end{equation}
In fact, by (\ref{3-24}) and (\ref{3-23}), we have
$$\ba{rl}
\ns&\ds\frac {1-e^{-2\hat\rho t}}{2\hat\rho}\|\Phi_{\widetilde{K}}(t,0)h\|^2=\ds\int^t_0 e^{-2\hat\rho (t-r)}\|\Phi_{\widetilde{K}}(t,0)h\|^2 dr\\
\ns&\le \ds\int^t_0 e^{-2\hat\rho (t-r)}\|\Phi_{\widetilde{K}}(t,r)\|^2\|\Phi_{\widetilde{K}}(r,0)h\|^2 dr\leq \ds\int^t_0 e^{-2\hat\rho(t- r)} \hat C_3^2e^{2\hat\rho(t-r)}\|\Phi_{\widetilde{K}}(r,0)h\|^2dr\\
\ns&<\hat C_3^2\ds\int^\infty_0\|\Phi_{\widetilde{K}}(r,0)h\|^2dr\leq \hat C_3^2\|Q^Z(0)\|\|h\|^2,\;\;\mbox{when}\;\; h\in H\;\;\mbox{and}\;\;t\geq 0,
\ea$$
 from which, it follows that
 $$
\ds\|\Phi_{\widetilde{K}}(t,0)\|\le \hat C_3\sqrt{\|Q^Z(0)\|{2\hat\rho}/(1-e^{-2\hat\rho t})}\;\;\mbox{for all}\;\; t>0.
$$
This, along with (\ref{WANGwgs3.76}), yields
$$\ba{rl}
\ns&\sup\limits_{t\in\mathbb{R}^+}\|\Phi_{\widetilde{K}}(t,0)\|=\max\big\{\max\limits_{t\in[0,T]}\|\Phi_{\widetilde{K}}(t,0)\|, \sup\limits_{t\in[T,\infty)}\|\Phi_{\widetilde{K}}(t,0)\|\big\}\\
\ns&\le \max\big\{  \max\limits_{t\in[0,T]}\|\Phi_{\widetilde{K}}(t,0)\|\,,~\hat C_3\sqrt{\|Q^Z(0)\|{2\hat\rho}/(1-e^{-2\hat\rho T})}\big\}=\hat C_4,
\ea$$
which leads to (\ref{wgs3.77-1}). Now, given $0\leq s\leq t$, there are only two possibilities: $(i)$ $[t/T]=[s/T]$; and $(ii)$ $[t/T]>[s/T]$.
In the first case,  (\ref{3-24-2}) follows from  (\ref{wgs3.73}) and (\ref{WANGwgs3.76}). In the second case, we have that $t\geq ([s/T]+1)T\geq s$.
This, along with (\ref{wgs3.64}), (\ref{wgs3.77-1}), (\ref{wgs3.73}) and (\ref{WANGwgs3.76}), indicates that
$$\ba{rl}\big\|\Phi_{\widetilde{K}}(t,s)\big\|&\le\big\|\Phi_{\widetilde{K}}\big(t,\big[{s}/{T}\big]T+T\big)\|\cdot\|\Phi_{\widetilde{K}}\big(\big[{s}/{T}\big]T+T,s\big)\big\|\\
\ns&=\big\|\Phi_{\widetilde{K}}\big(t-\big[{s}/{T}\big]T-T,0\big)\big\|\cdot\big\|\Phi_{\widetilde{K}}\big(T,s-\big[{s}/{T}\big]T\big)\big\|\leq \hat C_4 \hat C_2\leq \hat C_5.
\ea$$
In summary, we conclude that (\ref{3-24-2}) holds.\\

Finally, it follows from (\ref{3-24-2}) and (\ref{3-23}) that when $ t\geq 0$ and $ h\in H$,
$$t\big\|\Phi_{\widetilde{K}}(t,0)h\big\|^2=\ds\int^t_0\big\|\Phi_{\widetilde{K}}(t,0)h\big\|^2ds \le\ds\int^t_0\big\|\Phi_{\widetilde{K}}(t,s)\big\|^2\big\|\Phi_{\widetilde{K}}(s)h\big
\|^2ds\le \hat C_5^2\big\|Q^Z(0)\big\|\big\|h\big\|^2.$$
This implies that
$$\big\|\Phi_{\widetilde{K}}(t,0)\big\|\le \hat C_5\sqrt{{\|Q^Z(0)\|}/{t}},\;\;\mbox{when}\;\; t>0.$$
Therefore, there is an $N_0\in\mathbb{N}$ such that
$\|\Phi_{\widetilde{K}}(N_0T,0)\|<1$.
This, together with (\ref{wgs3.69}) and (\ref{wgs3.64}), yields that
$$
\hat\delta=\lim\limits_{k\rightarrow\infty}\big\|\Phi_{\widetilde{K}}(T,0)^{N_0k}\big\|^{\frac{1}{N_0k}}=\lim\limits_{k\rightarrow\infty}\big\|\Phi_{\widetilde{K}}(N_0T,0)^k\big\|^{\frac{1}{N_0k}}
\leq\big\|\Phi_{\widetilde{K}}(N_0T,0)\big\|^{\frac{1}{N_0}}<1.
$$
 This completes the proof. \endpf

\begin{proposition}\label{cor3.1}
When both $B(\cdot)$ and $D(\cdot)$ are time-invariant, i.e., $B(t)\equiv B$ and $D(t)\equiv D$ for all $t\geq 0$,  the following statements are equivalent:

\noindent $(i)$ Equation (\ref{state}) is linear $\hat T$-periodic feedback stabilizable for some $\hat T>0$.

\noindent $(ii)$ Equation (\ref{state}) is linear $\hat T$-periodic feedback stabilizable for any $\hat T>0$.

\noindent $(iii)$ Equation (\ref{state}) is linear time-invariant feedback stabilizable.
\end{proposition}

\noindent{\it Proof.} It suffices to show that $(i)\Rightarrow (iii)$. For this purpose, we suppose that $(i)$ holds.  Let $N\in\mathbb{N}$ with $N\ge 2$ and let
 $T=\hat T/N$. Since $B(\cdot)$ and $D(\cdot)$ are time-invariant, Equation (\ref{state}) is $T$-periodic. Because of $(i)$,
 there is an $NT$-periodic $\hat K^U_N(\cdot)\in L^\infty(\mathbb{R}^+;\mathcal{L}(H;U))$
such that the feedback equation (\ref{3-100}), where $Z=U$,
is exponentially stable.
  Now, by the same  way to show that Equation (\ref{state}) is LMPFS $\Rightarrow$ Equation (\ref{state}) is LPFS in the proof of Proposition \ref{lemma3.3} (where $Z=U$), we see that $\widetilde{K}(\cdot)$ given by
  (\ref{3-58-1}), where $Z=U$, is a LPFS law for
 Equation (\ref{state}). We claim that this $\widetilde{K}(\cdot)$ is time-invariant in the case that $B(\cdot)$ and $D(\cdot)$ are time-invariant. When this is done, $\widetilde{K}(\cdot)\equiv \widetilde{K}\in \cL(H; U)$  is a  feedback law for Equation (\ref{state}), which leads to $(iii)$.

 The rest is to show that  $\widetilde{K}(\cdot)$ is time-invariant. By the time-invariance of  $D(\cdot)$, and by (\ref{WUHENHEN5.54}) and  (\ref{3-58-1}), where $Z=U$, it suffices to show the valued function $W^U(t, h)$, given by (\ref{3.20}) with $Z=U$ is time-invariant. The later will be proved as follows.
Since Equation (\ref{state}) is time-invariant, we have that for each $t\in \mathbb{R}^+$, $h\in H$ and $u(\cdot)\in  L^2(\mathbb{R}^+; U)$,
$$y(s ;t,h,u)=y(s-t;0,h,v) \;\; \mbox{for all}\;\; s\geq t,$$
where $v(\cdot)$ is defined by $v(s)=u(s+t)$ for all $s\geq 0$. Hence, given  $t\in \mathbb{R}^+$ and $h\in H$,
$$\ds\int^\infty_t\left(\|y(s;t,h,u(s))\|^2 +\|u(s)\|_U^2\right)ds=\ds\int^\infty_0\left(\|y(r;0,h,u(r+t))\|^2 +\|u(r+t)\|_U^2\right)dr,$$
 for all  $u(\cdot)\in  L^2(\mathbb{R}^+; U)$.
Taking the infimum on the both sides of the above equation with respect to $u(\cd)\in L^2(\mathbb{R}^+; U)$, we get that $W^U(t,h)=W^U(0,h)$, i.e., the value function $W^U(t,h)$ is  time-invariant. This completes the proof.\endpf

\begin{remark}\label{remark3.1}
By Proposition~\ref{cor3.1}, we see that   linear time-periodic functions  $K(\cdot)$ will not aid in the linear stabilization of Equation (\ref{state}) when both $B(\cdot)$ and $D(\cdot)$ are time-invariant. On the other hand, when
 Equation (\ref{state}) is $T$-periodically time-varying,  linear time-periodic functions  $K(\cdot)$ do aid in the linear stabilization of Equation (\ref{state}). This can be seen from the following $2$-periodic ordinary differential equation:
$$y'(t)=\sum\limits_{j=1}^\infty\left[\chi_{[2j,2j+1)}(t)-\chi_{[2j+1,2j+2)}(t)\right]u(t).$$
For each $k\in\mathbb{R}$, consider the  feedback equation
$$y'(t)=\sum\limits_{j=1}^\infty\left[\chi_{[2j,2j+1)}(t)-\chi_{[2j+1,2j+2)}(t)\right]ky(t).$$
Clearly,  the corresponding Poincar$\acute{e}$ map  $\cP_k\equiv1$. Thus any linear time-invariant feedback equation
is not exponentially stable. On the other hand, by a direct computation, one can easily check that the 2-periodic time-varying feedback law given by
$$k(t)=\sum\limits_{j=1}^\infty\left[\chi_{[2j,2j+1)}(t)+2\chi_{[2j+1,2j+2)}(t)\right].$$
 is an LPFS law.

 \end{remark}

\section{The proof of Theorem~\ref{theorem1} and Theorem~\ref{theorem2}}

\noindent {\it Proof of Theorem \ref{theorem1}.}
$(a)\Leftrightarrow (b)$: We first show that $(b)\Rightarrow(a)$. Suppose that a subspace $Z$ of $U$ satisfies (\ref{3-14}). Then by Lemma \ref{lemma3.2},
  there is a finite dimension subspace $\hat Z$ of $Z $ such that
$ \hat V^{\hat Z}_{n_0}=H_1.$
From this, we can apply  Proposition~\ref{lemma3.1} to get  that
Equation (\ref{state}) is LMPFS with respect to $\hat Z$. This, along with  Proposition~\ref{lemma3.3}, yields that Equation (\ref{state}) is LPFS with respect to $\hat Z$.
Because $\hat Z$ is a subspace of $Z$,  Equation (\ref{state}) is  also LPFS with respect to $ Z$, i.e., $(a)$ stands.

We next show that $(a)\Rightarrow (b)$. Seeking for a contradiction, we  suppose that  $Z\in \mathcal{U}^{FS}$, but  (\ref{3-14}) does not hold. Then
  $\hat V^Z_{n_0}$  would be  a proper subspace of $H_1$.
 This, along with (\ref{2-17}), yields that $\hat V^Z$ is a proper subspace of $H_1$.
Write $(\hat V^Z)^\bot$ for the orthogonal complement subspace of $\hat V^Z$ in  $H$. Then, one can directly check that $H_1\bigcap(\hat V^Z)^\bot$ is the  orthogonal complement subspace of $\hat V^Z$ in  $H_1$, i.e.,
\begin{equation}\label{sgw3.70}
\big(H_1\bigcap (\hat V^Z)^\bot\big)\bot \hat V^Z;\;\;
H_1=\hat V^Z\bigoplus \big(H_1\bigcap (\hat V^Z)^\bot\big).
\end{equation}
Since $\hat V^Z$ is a proper subspace of $H_1$ and $\mbox{dim}H_1=n_0$ (see (\ref{wgs1.14})), we have
\begin{equation}\label{sgw3.71}
n_0\geq l\triangleq{\rm dim} \left(H_1\bigcap(\hat V^Z)^\bot\right)\geq 1.
\end{equation}
By (\ref{sgw3.70}) and (\ref{sgw3.71}),  we can let $\{\eta_1,\dots,\eta_{n_0}\}$ be a basis of $H_1$ such that $\{\eta_1,\cdots,\eta_l\}$  and $\{\eta_{l+1},\cdots,\eta_{n_0}\}$ are  bases of $H_1\bigcap(\hat
 V^Z)^\bot$   and  $ \hat V^Z $, respectively.
By (\ref{2-17}) in Lemma \ref{lemma2.3}, $\hat V^Z$ is an invariant subspace under $\cP_1$ .
 Thus
 there are matrices $A_1\in\mathbb{R}^{l\times l}$, $A_2\in\mathbb{R}^{(n_0-l)\times l}$, $A_3\in\mathbb{R}^{(n_0-l)\times (n_0-l)}$ such that
 \begin{equation}\label{3-15}
 \cP_1\left(\ba{cc} \eta_1, \cdots,\eta_l,\bigm|\eta_{l+1}, \cdots,\eta_{n_0}\ea\right)
 =\left(\ba{cc} \eta_1, \cdots,\eta_l,\bigm|\eta_{l+1},\cdots,\eta_{n_0}\ea\right)\left(\ba{cc} A_1 &0_{l\times(n_0-l)}\\A_2&
 A_3\ea\right).\end{equation}
 Let
  $P_{11}$ be the orthogonal  projection  from $H_1$ onto $H_1\bigcap(\hat V^Z)^\bot$.
 Define a linear bijection $\mathcal{J}: \mathbb{R}^l\rightarrow \big(H_1\bigcap(\hat V^Z)^\bot\big)$ by setting
 \begin{equation}\label{wwgs3.73}
 \mathcal{J}(\a)\deq(\eta_1,\dots,\eta_l)\a,\;\;\a\in \mathbb{R}^l,
 \end{equation}
 where $\a$ denotes the column vectors. By (\ref{3-15}) and (\ref{wwgs3.73}), we see that
   \begin{equation}\label{3-16}\ba{rl}
  \ns &P_{11}\cP_1^k \mathcal{J}(\a)=P_{11}\cP_1^k  (\eta_1,\cdots,\eta_l)\a\\
  \ns=&P_{11}(\eta_1,\cdots,\eta_l,\bigm|\eta_{l+1}, \cdots,\eta_{n_0})\left(\ba{cc} A_1 &0_{l\times(n_0-l)}\\A_2&
  A_3\ea\right)^k\left(\ba{l}\a\\0_{(n_0-l)\times 1}\ea\right)\\
  \ns=&(\eta_1,\cdots,\eta_l)A_1^k\a,\;\;\mbox{when}\;\; \a\in\mathbb{R}^l\;\;\mbox{and}\;\; k\in \mathbb{N}.
  \ea\end{equation}

 On the other hand, since $Z\in \mathcal{U}^{FS}$, there is a $T$-periodic $K(\cdot)\in L^\infty(\mathbb{R}^+;\cL(H;Z))$
   such that
 Equation (\ref{fstate}) is exponentially  stable, which implies that
 \begin{equation}\label{wgs3.75}
 \lim\limits_{t\rightarrow+\infty}y_K\left(t;0,h)\right)=0,\;\;\mbox{when}\;\; h\in H,
 \end{equation}
 where $y_K\left(\cdot;0,h)\right)$ denotes the solution of  Equation (\ref{fstate}) with the initial condition that $y(0)=h$.
 Let $u_K^h(t)=K(t)y_K(t; 0,h)$ for a.e. $t\geq 0$.
 The by Proposition 5.7 on Page 69 in \cite{Li}, we have
 \begin{equation}\label{3-17}
 y_K(t; 0,h)=\Phi(t,0)h+\ds\int^t_0\Phi(t,s)D(s) u^h_K(s)ds, \;\;\mbox{when}\;\;t\in\mathbb{R}^+\;\;\mbox{and}\;\;h\in H.
 \end{equation}
 From (\ref{3-17}) and (\ref{csubspace}), it follows that
\begin{equation}\label{wgs3.77}
P_{11}Py_K\left(kT;0,h)\right)\in P_{11}P (\cP^k h+ V^Z_k)\;\;\mbox{for all}\;\; h\in H\;\;\mbox{and}\;\; k\in \mathbb{N}.
\end{equation}
Since $P V^Z_k\triangleq \hat V^Z_k \subseteq\hat V^Z$ for all $k\in \mathbb{N}$ (see (\ref{2-14}) and (\ref{2-15})) and $P_{11}$ is the orthogonal projection from $H_1$ onto $H_1\bigcap(\hat V^Z)^\bot$, and because of (\ref{sgw3.70}), we have
$$
P_{11}P V^Z_k \subset P_{11}\hat V^Z=\{0\}.
$$
This, along with (\ref{wgs3.77}) and the fact that $P\cP^k=\cP^kP$ for all $k\in \mathbb{N}$ (see Parts $(a)$ and $(e)$ in Lemma~\ref{lemma2.2}), indicates that
\begin{equation}\label{wgs3.78}
P_{11}Py_K(kT;0,h)=P_{11}P\cP^kh=P_{11}\cP^kPh\;\;\mbox{for all}\;\; h\in H\;\;\mbox{and}\;\;k\in \mathbb{N}.
\end{equation}
Since $P: H\rightarrow H_1$ is a projection (see Lemma~\ref{lemma2.2}), it follows from (\ref{wgs3.75}) and (\ref{wgs3.78}) that
\begin{equation}\label{wgs3.79}
 \lim\limits_{k\rightarrow+\infty} P_{11}\cP^k h_1=0,\;\;\mbox{when}\;\;h\in H_1.
\end{equation}
Now by (\ref{wgs3.79}), (\ref{3-16}) and (\ref{2-13}), we have that $\lim\limits_{k\rightarrow\infty}A_1^k\a=0$, when $\a\in\mathbb{R}^l$,
from which, it follows that
\begin{equation}\label{wgs3.80}
\sigma(A_1)\in \mathbb{B}\;\;(\mbox{the open unit ball in}\;\; \mathbb{C}^1).
\end{equation}
By (\ref{3-15}), it holds that $\sigma(A_1)\subset\sigma(\cP_1)$. This, together with (\ref{wgs3.80}) and (\ref{WGSwgs2.22}), leads to a contradiction. Hence, $(a)\Rightarrow (b)$.
This completes the proof of $(a)\Leftrightarrow (b)$.

\vskip 5pt

$(b)\Leftrightarrow (c)$: First of all, we introduce two complex adjoint equations as follows:
\begin{equation}\label{a.1.4}
\psi'(t)-A^{*C}\psi(t)-B(t)^{*C}\psi(t) = 0\;\;\mbox{in}\;\; (0,n_0T),\;\; \psi(n_0T)\in H^C;
\end{equation}
\begin{equation}\label{a.1.5}
\psi'(t)-A^{*C}\psi(t)-B(t)^{*C}\psi(t) = 0\;\;\mbox{in}\;\; (0, T),\;\; \psi(T)\in H^C.
\end{equation}
 For each $\xi\in H^C$, Equation  (\ref{a.1.4}) (or (\ref{a.1.5})) with the initial condition that  $\psi^{\xi}_{n_0}(n_0T)={\xi}$ (or $\psi^{\xi}(T)={\xi}$) has a unique solution in $C[0,n_0T];H^C)$ (or $C([0,T]; H^C)$). We denote this solution by  $\psi^\xi_{n_0}(\cdot)$
 (or $\psi^\xi(\cdot)$). Clearly, when $\xi\in H$,   $\psi^\xi_{n_0}(\cdot)\in C[0,n_0T];H)$ and $\psi^\xi(\cdot)\in C([0,T]; H)$ are accordingly the solutions of (\ref{a.1.4}) and (\ref{a.1.5}) where $A^C$ and $B(t)^C$ are replaced by $A$ and $B(t)$ respectively.
 One can easily check that
\begin{equation}\label{a.1.51}
\psi^{\xi}(0)=\cP^{*C}\xi\;\; \mbox{and} \;\; \psi^{\xi}_{n_0}(0)=\left(\cP^{*C}\right)^{n_0}\xi \;\; \mbox{for all}\;\; \xi\in H^C.
\end{equation}
By the $T$-periodicity of $B^*(\cdot)$, we  see that for each $\xi\in H^C$,
\begin{equation}\label{a.1.7}
\psi^{\xi}_{n_0}((k-1)T+t)=\psi^{^{\xi_{k}}}(t),\; t\in[0,T],\,k\in\{1,\dots,n_0\},\;\;\mbox{where}\;\;\xi_{k}\deq\left(\cP^{*C}\right)^{n_0-k}\xi.
\end{equation}

Now we  carry out the proof of $(b)\Leftrightarrow (c)$ by several steps as follows.

\noindent {\it Step 1. To prove that  (\ref{3-14}) is  equivalent to the following condition:
\begin{equation}\label{a.1.8}
\forall\;  h\in H, \;\exists\;  u^h(\cd)\in L^2(\mathbb{R}^+ ;Z)\;\mbox{s.t.}\; Py(n_0T; 0,h, u^h)=0,\;\;\mbox{where}\;\; P\;\;\mbox{is given by}\;\;(\ref{wgs1.15}).
\end{equation}}
Suppose that (\ref{3-14}) holds. Then by (\ref{csubspace}), we have
\begin{equation}\label{a.1.9}
P\Big\{\ds\int^{n_0T}_0\Phi(n_0T,t)D(t)u(t)dt\bigm|u(\cd)\in L^2(\mathbb{R}^+ ;Z)\Big\}=H_1.
\end{equation}
Given $h\in H$, it holds that $P\Phi(n_0T,0)h\in H_1$ (see (\ref{2-2-1})). This, along with (\ref{a.1.9}), yields  that
there is a $u^h(\cd)\in L^2(\mathbb{R}^+ ;Z)$ such that
\begin{equation*}\label{a.1.10}
Py(n_0T;0,h, u^h)=P\Phi(n_0T,0)h+P\ds\int^{n_0T}_0\Phi(n_0T,t)D(t)u^h(t)dt=0,
\end{equation*}
which leads to  (\ref{a.1.8}).

Assume that  (\ref{a.1.8}) holds. Then  for any $h\in H$,  there exists $ u^h(\cd)\in L^2(\mathbb{R}^+ ;Z)$ such that $Py(n_0T;0,h, u^h)=0$, i.e.,
$$
-P\Phi(n_0T,0)h=P\ds\int_0^{n_0T}\Phi(n_0T,t)D(t)u^h(t)dt.
$$
Thus, we have
\begin{equation}\label{a.1.11}
\ba{rl}H_1&\supseteq \hat V_{n_0}^Z\deq P\Big\{\ds\int^{n_0T}_0\Phi(n_0T,t)D(t)u(t)dt\bigm|u(\cd)\in L^2(\mathbb{R}^+; Z)\Big\}\\
&\supseteq P\Big\{\ds\int^{n_0T}_0\Phi(n_0T,t)D(t)u^h(t)dt\bigm|h\in H\Big\}\\
&=-P\left\{\Phi(n_0T,0)h\bigm|h\in H\right\}=P\cP^{n_0}H.\ea
\end{equation}
By  the facts that $P\cP=\cP P$ (see (\ref{2-4})), $PH=H_1$  and $\cP H_1=\cP_1H_1=H_1$ (see (\ref{2-13}) and lemma \ref{lemma2.3}), we see that
$P\cP^{n_0}H=H_1$.
This, together with (\ref{a.1.11}), leads to (\ref{3-14}).

\medskip

\noindent {\it Step 2. To show that  $\xi\in P^*H_1$ and $\psi_{n_0}^\xi(0)=0$ $\Rightarrow$ $\xi=0$}

Recall  Lemma \ref{HUANGlemma2.3}. Because $\tilde H_1$ is a invariant subspace of $\cP^*$, it follows from
  (\ref{a.1.51}) that
\begin{equation}\label{a.1.52}
\psi_{n_0}^\xi(0)=(\cP^*)^{n_0}\xi=\big(\cP^*\bigl|_{\tilde H_1}\big)^{n_0}\xi\in\tilde H_1,\; \mbox{when}\; \xi\in \tilde H_1.
\end{equation}
By Lemma~\ref{HUANGlemma2.3}, we have
 $$
 \sigma\big(\cP^{*C}|_{\tilde H_1^C}\big)\bigcap\cB=\varnothing\;\;\mbox{and}\;\;\mbox{dim}\tilde H_1=n_0<\infty.
 $$
 Thus,  the map $\big(\cP^*\bigl|_{\tilde H_1}\big)^{n_0}$ is invertible from
$\tilde H_1$ onto $\tilde H_1$. Then by
(\ref{a.1.52}), we see that $\xi=0$ when $\xi\in\tilde H_1$ and
$\psi_{n_0}^\xi(0)=0$.
This, together with  (\ref{a.1.50}), implies that $\xi=0$ when $\xi\in P^* H_1$ and
$\psi_{n_0}^\xi(0)=0$.

\medskip

\noindent {\it Step 3. To show that (\ref{a.1.8}) $\Rightarrow$ (\ref{a.1.1})}

Clearly,  when $\eta,\,h\in H$ and $u(\cd)\in L^2(\mathbb{R}^+ ;Z)$,
\begin{equation}\label{a.1.12}
\big\langle \psi^\eta_{n_0}(0),\,h\big\rangle=\bigr\langle  \eta,\, y(n_0T; 0,h,u)\bigl\rangle
-\ds\int^{n_0T}_0\big\langle\bigr(D(t)\bigl|_Z\bigl)^*\psi^\eta_{n_0}(t), \,u(t)\big\rangle dt.
\end{equation}
Let $\xi\in P^*H_1$ satisfy the conditions on the left side of (\ref{a.1.1}).  Then by (\ref{a.1.12}) where $\eta=\xi$ and
$\psi^\xi_{n_0}(t)=\Phi(n_0T,t)^*\xi$, we find
\begin{equation}\label{a.1.13}
\big\langle \psi^\xi_{n_0}(0),\,h\big\rangle=\bigr\langle  \xi,\, y(n_0T; 0,h,u)\bigl\rangle,
{\mbox{ when }}h\in H{\mbox{ and }}u(\cd)\in L^2(\mathbb{R}^+; Z).
\end{equation}
By (\ref{a.1.8}), given $h\in H$, there is a $u^h(\cd)\in L^2(\mathbb{R}^+; Z)$ such that
\begin{equation}\label{a.1.53}
Py(n_0T;0,h, u^h)=0.
\end{equation}
Since $\xi\in P^*H_1$, there is  $g\in H_1$ such that $\xi=P^*g$.
This, along with
 (\ref{a.1.13}) and (\ref{a.1.53}), indicates
$$\big\langle \psi^\xi_{n_0}(0),\,h\big\rangle=\bigr\langle  \xi,\, y(n_0T; 0,h,u^h)\bigl\rangle
=\bigr\langle  P^*g,\, y(n_0T; 0,h,u^h)\bigl\rangle=\bigr\langle g,\, Py(n_0T; 0,h,u^h)\bigl\rangle=0,\;\forall\; h\in H.
$$
Hence,
$\psi^\xi_{n_0}(0)=0$. Then by the result in Step 2, we have $\xi=0$, i.e., (\ref{a.1.1}) holds.
\medskip

\noindent {\it Step 4. To show that (\ref{a.1.1}) $\Rightarrow$ (\ref{a.1.8})}\\
Assume that (\ref{a.1.1}) holds.
Define two subspaces
\begin{equation}\label{a.1.15}
\Gamma\deq \left\{\left(D(\cd)\bigl|_Z\right)^*\psi^{\xi}_{n_0} (\cd)\bigm|\xi\in P^*H_1\right\}\subseteq L^2(0,n_0T; Z)\; \mbox{and}\;
\Gamma_0\deq \left\{\psi^{\xi}_{n_0} (0)\bigm|\xi\in P^*H_1\right\}\subseteq H.
\end{equation}
Let
\begin{equation}\label{a.1.16}
\cL_1\left(\left(D(\cd)\bigl|_Z\right)^*\psi^{\xi}_{n_0} (\cd)\right)=\psi^{\xi}_{n_0} (0)
\;\;
{\mbox{ for all }}\;\;\xi\in P^*H_1.
\end{equation}
By (\ref{a.1.1}) and the result in Step 2, we see that $\cL_1:\Gamma\mapsto \Gamma_0$ is well defined. Clearly, it is linear.
 Given  $h\in H$, define  a linear functional $\cF^h$ on $\Gamma$ by
\begin{equation}\label{a.1.17}
\cF^h(\gamma)=\lan \cL_1(\gamma), h\ran\;\; {\mbox{ for all }}\;\;\gamma\in \Gamma.
\end{equation}
Since $\mbox{dim} (P^*H_1)=\mbox{dim} \tilde H_1=n_0<\infty$, it holds that $\mbox{dim} \Gamma<\infty$. Thus,   $\cF^h\in \cL(\Gamma; \mathbb{R})$.    By the Hahn-Banach theorem, there is a
$\bar \cF^h\in \cL(L^2(0,n_0T; Z); \mathbb{R})$
such that
\begin{equation}\label{a.1.18}
\bar \cF^h(\gamma)= \cF^h(\gamma)\;\; \mbox{for all}\;\;\gamma\in\Gamma;\;\;\mbox{and}\;\; \|\bar\cF^h\|=\|\cF^h\|.
\end{equation}
By making use of the Riesz representation theorem (see Page 59 in \cite{Diestel}), there
exists a function $u^h(\cdot)\in  L^2(0, n_0T;Z)$ such that
\begin{equation}\label{a.1.19}
\bar \cF^h(\gamma)=-\ds\int^{n_0T}_0\lan u^h(t), \gamma(t)\ran_Udt\;\;\mbox{for all}\;\; \gamma\in L^2(0, n_0T;Z).
\end{equation}
Now, since $P^*H_1=P^*H$ (see (\ref{a.1.50})), it follows from (\ref{a.1.16}), (\ref{a.1.17}), (\ref{a.1.18}) and (\ref{a.1.19}) that
\begin{equation*}\label{a.1.20}
-\ds\int^{n_0T}_0\big\langle\bigr(D(t)\bigl|_Z\bigl)^*\psi^{P^*\eta}_{n_0}(t), \,u^h(t)\big\rangle dt
=\big\langle \psi^{P^*\eta}_{n_0}(0),\,h\big\rangle \;\;\mbox{ for all }\;\;\eta\in H.
\end{equation*}
Meanwhile,  it follows by (\ref{a.1.12}) that for each $\eta\in H$,
\begin{equation*}\label{a.1.21}
\big\langle \psi^{P^*\eta}_{n_0}(0),\,h\big\rangle=\bigr\langle  {P^*\eta},\, y(n_0T; 0,h,u^h)\bigl\rangle
-\ds\int^{n_0T}_0\big\langle\bigr(D(t)\bigl|_Z\bigl)^*\psi^{P^*\eta}_{n_0}(t), \,u^h(t)\big\rangle dt
\end{equation*}
The above two equalities imply that
$$\bigr\langle  \eta,\, Py(n_0T; 0,h,u^h)\bigl\rangle=\bigr\langle  {P^*\eta},\, y(n_0T; 0,h,u^h)\bigl\rangle=0\;\;{\mbox{ for all }}\;\; \eta\in H,$$
i.e.,  $Py(n_0T; 0,h,u^h)=0$, which leads to  (\ref{a.1.8}).

From  Step 1-Step 4, one can easily check that  $(b)\Leftrightarrow (c)$.

\vskip 5pt

$(c)\Leftrightarrow (d)$: We first show that $(c)\Rightarrow (d)$. Suppose that a subspace $Z$ of $U$ satisfies  (\ref{a.1.1}).
Let $\mu$ and $\xi$ satisfy the conditions on the left side of  (\ref{a.1.2}) with the aforementioned $Z$. Then by (\ref{WANGJIAXU3.31}) (see Lemma~\ref{HUANGlemma2.3}),  it holds that $\xi\in {\tilde H_1^C}$. Hence, we can write
$\xi\deq \xi_1+i\xi_2$ with $\xi_1,\, \xi_2\in \tilde H_1$. By (\ref{a.1.50}), we have $\xi_1, \xi_2\in P^*H_1$.
By the last condition on the left side of  (\ref{a.1.2}), we have
\begin{equation*}\label{a.1.23}
  \left(D(t)\bigl|_Z\right)^{*C}\psi^{\xi}(t)=0\;\;\mbox{for a.e.}\;\; t\in(0,T).
\end{equation*}
By  (\ref{a.1.7}) and the third condition on the left side in (\ref{a.1.2}), it holds that
\begin{equation*}\label{a.1.24}
\psi^\xi_{n_0}((k-1)T+t)=\psi^{\mu^{n_0-k}\xi}(t)=\mu^{n_0-k}\psi^{\xi}(t) \;\;\mbox{for all}\;\; t\in [0,T], k=1,\dots, n_0.
\end{equation*}
Notice that $\psi^\xi_{n_0}(\cdot)=\psi^{\xi_1}_{n_0}(\cdot)+i\psi^{\xi_2}_{n_0}(\cdot)$. This, along with the above two equalities, yields that
$$\left(D(\cd)\bigl|_Z\right)^*\psi^{\xi_1}_{n_0}(\cd)+i\left(D(\cd)\bigl|_Z\right)^*\psi^{\xi_2}_{n_0}(\cd)=\left(D(\cd)\bigl|_Z\right)^*\psi^\xi_{n_0}(\cd)=0\;\;\mbox{over}\;\;(0, n_0T).$$
Since $\xi_1, \xi_2\in P^*H_1$, the above-mentioned equality, along with  (\ref{a.1.1}), leads to $\xi_1=\xi_2=0$, i.e., $\xi=0$. Hence, $Z$ satisfies  (\ref{a.1.2}). Thus, $(c)\Rightarrow (d)$.

We next show that $(d)\Rightarrow (c)$.  Suppose that a subspace $Z$ satisfies (\ref{a.1.2}).
 In order to show that $Z$ satisfies (\ref{a.1.1}), it suffices to prove
\begin{equation}\label{a.1.27}
 \hat\xi\in (P^*H_1)^C\;\;\mbox{and}\;\; \left(D(\cd)\bigl|_Z\right)^{*C}\psi_{n_0}^{\hat \xi}(\cd)=0\;\;\mbox{over}\;\; (0, n_0T)\Rightarrow \hat\xi=0.
\end{equation}
First of all, we notice   that $(P^*H_1)^C=\tilde H_1^C$ and $\mbox{dim}\tilde H_1^C=n_0$  (see Lemma~\ref{HUANGlemma2.3}).
In this step, we simply write
$$\cQ\deq \cP^{*C}\bigl|_{\tilde H_1^C}\in\cL(\tilde H_1^C)\;\;\mbox{and}\;\;
 D_1(\cd)\deq \big(\big(D(\cd)\bigl|_Z\big)^{*C}\big)\bigl|_{(0,T)}\in L^2(0,T ;\cL(H,Z)).$$
By Lemma~\ref{HUANGlemma2.3} and (\ref{a.1.42}), we have that  $\sigma (\cQ) =\{\bar \l_j\}_{j=1}^n$;
 $l_j$ is the algebraic multiplicity of $\bar\l_j$.
Hence, $\mathbb{P}_0(\l)\triangleq\prod\limits^{n}_{j=1}\left(\l-\bar\l_j\right)^{l_j}$ is the characteristic polynomial of $\cQ$.
Write $\hat l_j$ for the geometric multiplicity of $\bar\l_j$.
Clearly,
$\hat l_j\le l_j$ for all $j$.
Let $\beta\triangleq\{\beta_1,\dots,\beta_{n_0}\}$ be a basis of $(P^*H)^C=\tilde H_1^C$ such that
\begin{equation}\label{a.1.25}
\cQ(\beta_1,\dots,\beta_{n_0})=J(\beta_1,\dots,\beta_{n_0}).
\end{equation}
Here $J$ is the Jordan matrix:
 ${\rm diag }\big\{J_{11}, \dots,J_{1\hat l_1},J_{21}, ~\dots,~J_{2\hat l_2},\dots,J_{n1}, ~\dots,~J_{n\hat l_n}\big\}$ with
 $$J_{jk}=\left(\ba{cccc}
 \bar\l_j&1\\&\ddots&\ddots\\&&\ddots&1\\&&&\bar\l_j
 \ea\right) \;\;\mbox{a}\;\; d_{jk}\times d_{jk}\;\;\mbox{matrix},$$
 where $j=1,\dots, n$,  $k=1,\dots, \hat l_j$, and for each $j$, $\{d_{jk}\}_{k=1}^{\hat l_j}$ is decreasing. It is clear that
 $\sum\limits^{\hat l_j}_{k=1}d_{jk}= l_j$ for each $ j=1,\dots, n$, and $\sum\limits^{n}_{j=1}\sum\limits^{\hat l_j}_{k=1}d_{jk}=n_0$. We rewrite the basis $\beta$ as
 $$\beta=\big\{\xi_{111},\dots,\xi_{11d_{11}},\xi_{1\hat l_11},\dots,\xi_{1\hat l_1d_{1\hat l_1}} ,\dots,
 \xi_{n11},\dots,\xi_{n1d_{n1}},\xi_{n\hat l_n1},\dots,\xi_{n\hat l_nd_{n\hat l_n}} \big\}.$$
Then by (\ref{a.1.25}), one can easily check that for each $j\in\{1,\dots,n\}$ and $k\in\{1,\dots, \hat l_j\}$,
 \begin{equation}\label{a.1.26}
 \ba{r}
 \ns\big(\bar\l_jI-\cQ\big)^q\xi_{jkr}=\left\{
 \ba{cl}
 \xi_{jk(r-q)} \q & {\rm  when}~r>q,\\ 0 &{\rm  when}~r\le q.
 \ea
 \right.\ea
 \end{equation}

Now we assume $\hat\xi$ satisfies the conditions on the left side of (\ref{a.1.27}). Since $\hat\xi\in (P^*H_1)^C=\tilde H_1^C$, there is a vector
$$ \big(C_{111},\dots,C_{11d_{11}},C_{1\hat l_11},\dots,C_{1\hat l_1d_{1\hat l_1}} ,\dots,
C_{n11},\dots,C_{n1d_{n1}},C_{n\hat l_n1},\dots,C_{n\hat l_nd_{n\hat l_n}} \big)^*\in \mathbb{C}^{n_0},$$
such that
\begin{equation}\label{a.1.28}
\hat \xi=\sum\limits^n_{j=1}\sum\limits^{\hat l_j}_{k=1}\sum\limits^{d_{jk}}_{r=1}C_{jkr}\xi_{jkr}.
\end{equation}
From (\ref{a.1.51}) and the second condition on the left side of (\ref{a.1.27}), it follows that for each $m\in\{0,\dots, n_0-1\}$,
 $$ D_1(\cd)\psi_{n_0}^{\hat \xi}(\cd)\bigl|_{((n_0-m-1)T,(n_0-m)T)}=0\;\;\mbox{i.\,e.}\;\;\sum\limits^n_{j=1}\sum\limits^{\hat l_j}_{k=1}\sum\limits^{d_{jk}}_{r=1}C_{jkr} D_1(t)
\psi^{\cQ^m\xi_{jkr}}(t)=0\;\;\mbox{for a.e.}\;\;t\in (0,T),  $$
from which,  we see that
\begin{equation}\label{a.1.29}\ba{c}
\ns\sum\limits^n_{j=1}\sum\limits^{\hat l_j}_{k=1}\sum\limits^{d_{jk}}_{r=1}C_{jkr} D_1(\cd)
\psi^{g(\cQ)\xi_{jkr}}=0\;\;\mbox{over}\;\;(0,T)
\ea\end{equation}
for any polynomial $g$ with $\mbox{degree} (g)\leq n_0-1$.
Arbitrarily fix a  $\widetilde{j}\in\{1,\dots, n\}$. Let
$$\mathbb{P}_{\widetilde{j}}(\l)=\prod\limits _{j=1, j\neq \widetilde{j}}^n\left(\l-\bar\l_j\right)^{l_j}.
$$
By taking $g(\l)=\l^m \mathbb{P}_{\widetilde{j}}(\l)$, with $m=0,1,\dots, l_{\tilde j}-1$, in (\ref{a.1.29}),  we have
\begin{equation*}\label{a.1.30}
\sum\limits^n_{j=1}\sum\limits^{\hat l_j}_{k=1}\sum\limits^{d_{jk}}_{r=1}C_{jkr} D_1(\cd)
\psi^{\cQ^m\mathbb{P}_{\widetilde{j}}(\cQ)\xi_{jkr}}(\cd)=0\;\;\mbox{over}\;\; (0,T),\;\;\mbox{when}\;\; m\in\{0,1,\dots, l_{\widetilde{j}}-1\}.
\end{equation*}
By (\ref{a.1.26}), we see that
$$\mathbb{P}_{\widetilde{j}}(\cQ)\xi_{jkr}=0,\;\;\mbox{ when }\;\;j\in\{1,\dots,,n\}, j\neq\widetilde{j},\,k\in\{1,\dots, \hat l_j\},\,r\in\{1,\dots, d_{jk}\}.$$
The above two equalities imply that for each $m\in\{0,1,\dots, l_{\widetilde{j}}-1\}$,
$$\sum\limits^{\hat l_{\widetilde{j}}}_{k=1}\sum\limits^{d_{{\widetilde{j}}k}}_{r=1}C_{{\widetilde{j}}kr} D_1(\cd)
\psi^{\cQ^m\mathbb{P}_{\widetilde{j}}(\cQ)\xi_{{\widetilde{j}}kr}}(\cd)=0\;\;\mbox{over}\;\;(0,T),
$$
from which, it follows that
\begin{equation}\label{a.1.31}\ba{c}
\ns\sum\limits^{\hat l_{\widetilde{j}}}_{k=1}\sum\limits^{d_{{\widetilde{j}}k}}_{r=1}C_{{\widetilde{j}}kr} D_1(\cd)
\psi^{f(\cQ)\mathbb{P}_{\widetilde{j}}(\cQ)\xi_{{\widetilde{j}}kr}}(\cd)=0\;\;\mbox{over}\;\;(0,T),
\ea\end{equation}
for any polynomial $f$ with $\mbox{degree} (f)\leq l_{\tilde j}-1$.
Given $m\in\{0,1,2,\cdots, l_{\widetilde{j}}-1\}$, since $\mathbb{P}_{\widetilde{j}}(\l)$ and $(\l-\bar\l_{\widetilde{j}})^{m+1}$ are  coprime,
 there are  polynomials $g^1_m (\l)$ and $g^2_m(\l)$ with $\mbox{degree} (g^1_m)\leq m$ and $\mbox{degree} (g^2_m)\leq \mbox{degree} (\mathbb{P}_{\widetilde{j}})-1$, respectively,  such that
 $$
 g^1_m(\l)\mathbb{P}_{\widetilde{j}}(\l)+g^2_m(\l)(\l-\bar\l_{\widetilde{j}})^{m+1} \equiv 1,
 $$
from which, we see that
 \begin{equation}\label{a.1.32}
(\cQ-\bar\l_{\widetilde{j}}I)^{l_{\widetilde{j}}-m-1} g^1_m(\cQ)\mathbb{P}_{\widetilde{j}}(\cQ)\xi_{{\widetilde{j}}kr}\\
+g^2_m(\cQ)(\cQ-\bar\l_{\widetilde{j}}I)^{l_{\widetilde{j}}}\xi_{{\widetilde{j}}kr}
 \equiv (\cQ-\bar\l_{\widetilde{j}}I)^{l_{\widetilde{j}}-m-1}\xi_{{\widetilde{j}}kr},
 \end{equation}
for all $m\in\{0,1,\cdots, l_{\widetilde{j}}-1\}$, $k\in\{1,2,\cdots,\hat l_{\widetilde{j}}\}$, and $r\in\{1,2,\cdots,d_{\widetilde{j}k}\}$.
By (\ref{a.1.26}), we have
\begin{equation}\label{wang5.37}
(\hat\cQ-\bar\l_{\widetilde{j}}I)^{l_{\widetilde{j}}}\xi_{{\widetilde{j}}kr}=0\;\;\mbox{for all}\;\; k\in\{1,2,\cdots, \hat l_{\widetilde{j}}\},\,r\in\{1,2,\cdots, d_{{\widetilde{j}}k}\}.
\end{equation}
Taking  $f(\l)=(\l-\bar\l_{\widetilde{j}})^{l_{\widetilde{j}}-m-1}g^1_m(\l)$, with  $m=0,\dots,l_{\widetilde{j}}-1$, in (\ref{a.1.31}),  using (\ref{a.1.32}) and (\ref{wang5.37}), we find

\begin{equation}\label{a.1.33}
\sum\limits^{\hat l_{\widetilde{j}}}_{k=1}\sum\limits^{d_{{\widetilde{j}}k}}_{r=1}C_{{\widetilde{j}}kr} D_1(\cd)
\psi^{(\cQ-\bar\l_{\widetilde{j}}I)^m\xi_{{\widetilde{j}}kr}}(\cd)=0\;\;\mbox{over}\;\;(0,T),\;\;\mbox{for each}\;\; m\in\{0,1,\dots, l_{\widetilde{j}}-1\}.
\end{equation}

Now we are on the position  to show
\begin{equation}\label{a.1.55}
C_{\widetilde{j}kr}=0\;\;\mbox{for all}\;\;k\in\{1,2,\cdots, \hat l_{\widetilde{j}}\},\,r\in\{1,\dots, d_{{\widetilde{j}}k}\},
\end{equation}
which leads to $\hat\xi=0$ because of (\ref{a.1.28}). For this purpose, we write
$$K_{\widetilde{j}}^m=\left\{k\in\{1,2,\dots, \hat l_{\widetilde{j}}\}\bigm|
d_{{\widetilde{j}}k}> m\right\},\;\; m=0,1,\dots, l_{\widetilde{j}}-1.
$$
One can easily check that (\ref{a.1.55}) is equivalent to
\begin{equation}\label{a.1.35}
\mathcal{C}_{\hat m}\deq\left\{C_{\widetilde{j}k\hat m},~~k\in K_{\widetilde{j}}^{\hat m-1}\right\}=\{0\}\;\;\mbox{for all}\;\;
\hat m\in\{1,\dots, d_{\widetilde{j}1}\}
\end{equation}
We will use the  mathematical induction method with respect to $\hat m$  to prove (\ref{a.1.35}). (Notice that
$d_{\widetilde{j}k}$ is decreasing with respect to $k$.)  First of all, we let
\begin{equation}\label{a.1.34}
\mathbb{Q}_{\widetilde{j}}^m(\l)=\left(\bar\l_{\widetilde{j}}-\l\right)^{m},\;\;
 m=0,1,\dots, l_{\widetilde{j}}-1,
\end{equation}
 In the case that  $\hat m=d_{\widetilde{j}1}$, it follows from  (\ref{a.1.34}) and (\ref{a.1.26}) that
$$
\mathbb{Q}^{\hat m-1}_{\widetilde{j}}(\cQ)\xi_{\widetilde{j}k\hat m}=\big(\bar\l_{\widetilde{j}}I-\cQ\big)^{\hat m-1}\xi_{\widetilde{j}k\hat m}=\xi_{\widetilde{j}k1},\;\;\mbox{when}\;\; k\in K^{\hat m-1}_{\widetilde{j}},\emph{}$$
and
$$
\mathbb{Q}^{\hat m-1}_{\widetilde{j}}(\cQ)\xi_{\widetilde{j}kr}=0,\q {\mbox{ when }}
 k\in K^{\hat m-1}_{\widetilde{j}},\, r<\hat m; {\mbox{ or }}k\notin K^{\hat m-1}_{\widetilde{j}}, \, r\in\{1,\dots, d_{\widetilde{j}k}\}.$$
These, alone with (\ref{a.1.33}) (where  $m=\hat m-1$),  imply that
\begin{equation*}\label{a.1.36.1}
\sum\limits_{k\in K_{\widetilde{j}}^{\hat m-1}}C_{\widetilde{j}k\hat m} D_1(\cd)
\psi^{\xi_{{\widetilde{j}}k1}}(\cd)=0\;\;\mbox{over}\;\;(0,T).
\end{equation*}
Let
$$
\bar\xi_{\hat m}\deq\sum_{k\in K_{\widetilde{j}}^{\hat m-1}}C_{{\widetilde{j}}k\hat m}\xi_{{\widetilde{j}}k1},\; \hat m=1,\dots, d_{\widetilde{j}1}.
$$
Then, it holds that
\begin{equation}\label{wang5.42}
 D_1(\cd)
\ds\psi^{\bar\xi_{\hat m}}(\cd)=0\;\;\mbox{over}\;\;(0,T).\end{equation}
Since for each $k\in\{1,\dots,\hat l_{\tilde j}\}$,
$\xi_{\widetilde{j}k1}$ is an eigenfunction of $\cQ$ with respect to the eigenvalue $\bar\l_{\widetilde{j}}$, it follows from the definition of $\bar\xi_{\hat m}$ that $(\bar\l_{\widetilde{j}}I-\cQ)\bar\xi_{\hat m}=0$. This, along with (\ref{wang5.42}) and (\ref{a.1.2}),
yields that $\bar\xi_{\hat m}=0$, i.e.,
 $\bar\xi_{d_{\widetilde{j},1}}=0$, which leads to $\mathcal{C}_{d_{\tilde j1}}=0$ because of the
 linear independence of the group $\{\xi_{{\widetilde{j}}k1},~k\in K_{\widetilde{j}}^{\hat m-1}\}$. Hence,  (\ref{a.1.35}) holds when $\hat m=d_{\widetilde{j}1}$.

\medskip

Suppose inductively that  (\ref{a.1.35}) holds when  $ \tilde m+1\le\hat m\le d_{\widetilde{j}1}$ for some $\tilde m\in\{1,\dots,d_{\widetilde{j}1}-1\}$,
 i.e.,
 \begin{equation}\label{wang5.43}
 \mathcal{C}_{\hat m}=\{0\},\;\;\mbox{when}\;\;\tilde m+1\le\hat m\le d_{\widetilde{j}1}.
 \end{equation}
   We will show that
 (\ref{a.1.35}) holds when  $\hat m=\tilde m$, i.e., $\mathcal{C}_{\tilde m}=\{0\}$.
 In fact, it follows from (\ref{a.1.26}) that
$$
\mathbb{Q}^{\tilde m-1}_{\widetilde{j}}(\cQ)\xi_{\widetilde{j}kr}=\left\{\ba{cl}
\ns\xi_{\widetilde{j}k(r-\tilde m+1)},\q &{\rm when }~ k\in K^{\tilde m-1}_{\widetilde{j}}, \, r\ge\tilde m,\\
\ns0,\q &{\rm when }~ k\in K^{\tilde m-1}_{\widetilde{j}}, \, r<\tilde m,\\
\ns0,\q &{\rm when }~ k\notin K^{\tilde m-1}_{\widetilde{j}}, \, r\in\{1,\dots, d_{\widetilde{j}k}\}.
\ea\right.
$$
 This, alone with (\ref{a.1.33}) (where  $m=\tilde m-1$),  indicates that
$$\sum\limits^{\hat l_{\widetilde{j}}}_{k=1}\sum\limits^{d_{{\widetilde{j}}k}}_{r=1}C_{{\widetilde{j}}kr} D_1(\cd)
\psi^{\mathbb{Q}^{\tilde m-1}_{\widetilde{j}}(\cQ)\xi_{{\widetilde{j}}kr}}(\cd)
=\sum\limits_{k\in K_{\widetilde{j}}^{\tilde m-1}}\sum\limits^{d_{{\widetilde{j}}k}}_{r=\tilde m}C_{{\widetilde{j}}kr} D_1(\cd)
\psi^{\xi_{{\widetilde{j}}k(r-\tilde m+1)}}(\cd)=0\;\;\mbox{over}\;\;(0,T).$$
Then, by (\ref{wang5.43}), we have
\begin{equation}\label{a.1.36}
\sum\limits_{k\in K_{\widetilde{j}}^{\tilde m-1}}C_{\widetilde{j}k\tilde m} D_1(\cd)
\psi^{\xi_{{\widetilde{j}}k1}}(\cd)=0\;\;\mbox{over}\;\;(0,T).\end{equation}
Let
$$
\bar\xi_{\tilde m}\deq\sum_{k\in K_{\widetilde{j}}^{\tilde m-1}}C_{{\widetilde{j}}k\tilde m}\xi_{{\widetilde{j}}k1}.
$$
Then, it follows from (\ref{a.1.36}) that
\begin{equation}\label{a.1.37}
 D_1(\cd)
\ds\psi^{\bar\xi_{\tilde m}}(\cd)=0\;\;\mbox{over}\;\;(0,T).\end{equation}
Since for each $k\in\{1,\dots,\hat l_{\widetilde{j}}\}$, $\xi_{\widetilde{j},k,1}$ is an eigenfunction of $\cQ$ with respect to the eigenvalue $\bar\l_{\widetilde{j}}$, it holds that $(\bar\l_{\widetilde{j}}I-\cQ)\bar\xi_{\tilde m}=0$ . This, along with (\ref{a.1.37}) and (\ref{a.1.2}), yields that
  $\bar\xi_{\hat m}=0$. Hence, $\mathcal{C}_{\tilde m}=\{0\}$ because of  the linear independence of the group
$\{\xi_{{\widetilde{j}}k1},~k\in K_{\widetilde{j}}^{\tilde m-1}\}$.
In summary, we conclude that $(d)\Rightarrow (c)$.
 This completes the proof of Theorem \ref{theorem1}. \endpf

\bigskip

\noindent{\it Proof of Theorem \ref{theorem2}.} Clearly, it suffice to show the {\it only if} part. Assume that Equation (\ref{state}) is LPFS. By the equivalence of $(a)$ and $(b)$ in Theorem \ref{theorem1}, it holds that
\begin{equation}\label{wgs3.81}
\hat V^U_{n_0}=H_1.
\end{equation}
Meanwhile, according to  Lemma \ref{lemma3.2}, there is a finite dimensional subspace of $\hat Z$ of  $U$, such that
\begin{equation}\label{wgs3.82}
\hat V^U_{n_0}=\hat V^{\hat Z}_{n_0}.
\end{equation}
  From (\ref{wgs3.81}) and (\ref{wgs3.82}), it follows that
    $\hat V^{\hat Z}_{n_0}=H_1$.
 This, along with the equivalence of $(a)$ and $(b)$ in Theorem \ref{theorem1}, indicates that Equation (\ref{state})  is LPFS with respect to $\hat Z$.
 Hence, we complete the proof of Theorem \ref{theorem2}.
 \endpf

 \bigskip

\section{Applications}
In this  section, we will present some applications of Theorem \ref{theorem1} to the internally controlled heat equations with time-periodic potential.

Let  $\Omega$ be a bounded domain in $\mathbb{R}^d$ ($d\geq 1$) with a $C^2$-smooth boundary $\partial\Omega$. Write $Q\triangleq\Omega\times\mathbb{R}^+$ and $\sum\triangleq\partial \Omega\times \mathbb{R}^+$. Let  $\omega\subseteq \Omega$  be a non-empty open subset with its characteristic function $\chi_{\omega}$.
Let $T>0$ and   $a\in L^\infty (Q)$ be   $T$-periodic (with respect to the time variable  $t$), i.e., for a.e. $t\in \mathbb{R}^+$, $a(\cdot,t)=a(\cdot,t+T)$ over $\Omega$.
 One can easily check that the function $a$ can be treated as a $T$-periodic function in    $ L^1_{loc}(\mathbb{R}^+;L^2(\Omega)).$
   Consider the following controlled heat equation:
\begin{equation}\label{hstate} \left\{\begin{array}{lll}
\ns \partial_ty(x,t)-\triangle y(x,t)+a(x,t)y(x,t)=\chi_\omega(x)u(x,t)&\mbox{in}&Q,\\
\ns y(x,t)=0&\mbox{on}&\sum,
\end{array}\right.\end{equation}
where $u\in L^2(\mathbb{R}^+;L^2(\Omega))$.
Given $y_0\in L^2(\Omega)$ and  $u\in L^2(\mathbb{R}^+;L^2(\Omega))$, Equation (\ref{hstate}) with the initial condition that $y(x,0)=y_0(x)$ has a unique solution
$y(\cdot; 0,y_0,u)\in C(\mathbb{R}^+; L^2(\Omega))$.

 Let $H=U=L^2(\Omega)$ and $A=-\triangle$ with $\mathcal{D}(A)=H_0^1(\Omega)\bigcap H_2(\Omega)$. Define, for a.e. $t\in \mathbb{R}^+$, $B(t): H \rightarrow H$ by $B(t)z(x)=a(x,t)z(x)$, $x\in\Omega$,
 and $D(t): U\rightarrow H$ by $D(t)v(x)=\chi_\omega(x)v(x)$, $x\in\Omega$. Clearly, $(-A)$ generates a compact semigroup on $L^2(\Omega)$ and both $B(\cdot)\in L^1_{loc}(\mathbb{R}^+;L^2(\Omega))$
 and $D(\cdot)\in L^\infty(\mathbb{R}^+; \cL(U;H))$ are $T$-periodic. Thus, we can study Equation (\ref{hstate}) under the framework (\ref{state}).
 Write
  $\{\Psi_a(t,s)\}_{0\leq s\leq t}$ for the evolution system generated by $-A-B(\cdot)$. We  use   notations $n_0$, $P$, $H_j$ (with $j=1,2$),  $V^Z_k$ and $\hat V^Z_k$ (with $k\in \mathbb{N}$) to denote the same subjects
as those introduced in section 1. We will use the different equivalent conditions in Theorem \ref{theorem1} to show that Equation (\ref{hstate}) is LPFS.

\begin{corollary} Equation (\ref{hstate}) is LPFS. Consequently, it is LPFS with respect to a finite dimensional subspace of $L^2(\Omega)$.

\end{corollary}
{\it Proof.} We will provide two ways to show that Equation (\ref{hstate}) is LPFS.
 We first use the equivalence $(a)\Leftrightarrow(c)$ in Theorem \ref{theorem1}. In fact,
$\psi(\cdot)\triangleq \Psi_a(n_0T, \cdot)^*\xi$ with $\xi\in H$ is the solution to the equation:
\begin{equation}\label{a.2.7}
\left\{\begin{array}{lll}
\ns \partial_t\psi(x,t)+\triangle \psi(x,t)-a(x,t)\psi(x,t)=0&\mbox{in}&\Omega\times (0,n_0T),\\
\ns \psi(x,t)=0&\mbox{on}&\partial\Omega\times(0,n_0T),\\
\ns \psi(x,n_0T)=\xi(x)&\mbox{in}&\Omega,
\end{array}\right.\end{equation}
and it holds that
\begin{equation}\label{a.2.3}
D(t)\eta=\chi_{\omega}\eta\;\;\mbox{for any}\;\;\eta\in H\;\;\mbox{and}\;\; t\in[0,T].
\end{equation}
These, along with the unique continuation  property of parabolic equations established in \cite{LFH} (see also \cite{Phung1} and \cite{Phung2}), leads to  the condition $(c)$ in Theorem \ref{theorem1} for the current case. Then, according to
 the equivalence $(a)\Leftrightarrow(c)$ in Theorem \ref{theorem1}, Equation (\ref{hstate}) is LPFS.

We next use the equivalence $(a)\Leftrightarrow (b)$ in Theorem \ref{theorem1}.
 Without loss of generality, we can assume that $n_0\geq 1$, for otherwise Equation (\ref{hstate}), with the null control $u=0$, is stable. When $n_0\geq 1$, we have $H_1\neq\{0\}$ and $\|P\|>0$.
Write
$\left\{\xi_1,\dots,\xi_{n_0}\right\}$ for an orthonormal basis of $H_1$.
By the  approximate controllability of the heat equation (\see \cite{Fabra}),
$V^U_1$ is dense in $H$. Thus there are  $\eta_j$, $j=1\dots, n_0$, in $ V^U_1$ such that
such that
\begin{equation}\label{4-1-1}
\|\eta_j-\xi_j\|\le\displaystyle\frac{1}{16n_0\|P\|}\;\;\mbox{for all}\;\; j=1,\dots, n_0.
\end{equation}
Since $P$ is a projection from $H$ onto $H_1$, we have $P\xi_j=\xi_j$ for all $j=1,\dots,n_0$. This, along with (\ref{4-1-1}), yields that for each $j\in \{1,\dots, n_0\}$,
\begin{equation}\label{4-2-1}\ba{rl}
\ns\|P\eta_j\|\le &\|\xi_j\|+\|P\eta_j-\xi_j\|=\|\xi_j\|+\|P(\eta_j-\xi_j)\|\\
\ns\le&\|\xi_j\|+\|P\|\|\eta_j-\xi_j\|\le1+\ds{1\over{16n_0}};
\ea\end{equation}
and
\begin{equation}\label{4-2-2}\ba{rl}
\ns\lan P\eta_j,\xi_j\ran=&\lan\xi_j+(P\eta_j-\xi_j),\xi_j\ran=1+\lan P\eta_j-\xi_j,\xi_j\ran\\
\ns\ge& 1-\|P\eta_j-\xi_j\|\ge 1-\|P\|\|\eta_j-\xi_j\|\ge1-\ds\frac{1}{16n_0}.\\
\ea\end{equation}
Since  $P\eta_j\in H_1$ and $\{\xi_k\}_{k=1}^{n_0}$ is an orthonormal basis of $H_1$, it holds that
\begin{equation}\label{4-2-3}
\|P\eta_j\|^2=\sum^{n_0}_{k=1}|\lan P\eta_j,\xi_k\ran|^2,\;\;\mbox{when}\;\;j=1,\dots,n_0.
\end{equation}
From (\ref{4-2-3}), (\ref{4-2-1})  and (\ref{4-2-2}), it follows that for each $j\in \{1,\dots, n_0\}$,
$$
\ba{rl}
\sum\limits_{k\neq j}\left|\lan P\eta_j,\xi_k\ran\right|&\le(n_0-1)^{1/2} \left(\sum\limits_{k\neq j}\left|\lan P\eta_j,\xi_k\ran\right|^2\right)^{1/2}\\
\ns &=(n_0-1)^{1/2}\big(\|P\eta_j\|^2-|\lan P\eta_j,\xi_j\ran|^2\big)^{1/2}\\
\ns&\le
n_0^{1/2}\big(\left(1+{1}/(16n_0)\right)^2-\left(1-{1}/(16n_0)\right)^2\big)^{1/2}={1}/{2}.
\ea$$
This, together with (\ref{4-2-2}), indicates that
\begin{equation}\label{4-3}
\lan P\eta_j,\xi_j\ran\ge 1-{1}/(16n_0)>{1}/{2}\ge \sum_{k\neq j}\left|\lan P\eta_j,\xi_k\ran\right|,\; j=1,2,\cdots, n_0.
\end{equation}

We claim that $\{P\eta_1,\dots,P\eta_{n_0}\}$ is a  linearly independent group.
In fact, suppose that
\begin{equation}\label{WGS4.8}
\sum\limits_{j=1}^n c_j P\eta_j=0\;\;\mbox{for some}\;\; c_1,\dots, c_{n_0}\in\mathbb{R}.
\end{equation}
Write $\hat A\triangleq(\lan P\eta_j,\xi_k\ran)_{j,k}\in \mathbb{R}^{n_0\times n_0}$ and  $\hat c\triangleq(c_1,\dots, c_{n_0})^*\in\mathbb{R}^{n_0}$.
By (\ref{4-3}), the matrix $\hat A$ is diagonally dominant, hence it is invertible. Then, from (\ref{WGS4.8}), it follows that
$ \hat A^*\hat c=0,$ which implies $\hat c=0$. Hence, $P\eta_1,\dots,P\eta_{n_0}$ are linearly independent.

Since $\mbox{dim} H_1=n_0$, it follows that
$\mbox{span}\{P\eta_1,\dots,P\eta_{n_0}\}=H_1$.
Therefore
$$H_1\supseteq\hat V^U_{n_0}\supseteq\hat V^U_1=PV^U_1\supseteq \mbox{span}\{P\eta_1,\cdots,P\eta_{n_0}\}=H_1.
$$
from which, it follows that $H_1=\hat V^U_{n_0}$. This, along with the equivalence of $(a)$ and $(b)$ in Theorem \ref{theorem1}, indicates that Equation (\ref{hstate}) is LPFS.

Finally,  according to  Theorem \ref{theorem2}, there is a finite-dimensional subspace $Z$ of $U$ such that
Equation (\ref{hstate}) is LPFS with respect to $Z$. This completes the proof.

\bigskip

\begin{corollary} Equation (\ref{hstate}) is LPFS
with respect to the subspace $P^*H$.
\end{corollary}

\noindent{\it Proof.}
Let $Z=P^*H$.  By the equivalence between $(a)$ and $(d)$ in Theorem~\ref{theorem1}, it suffices to show that $Z$ satisfies (\ref{a.1.2}), i.e.,
\begin{equation}\label{a.2.1}
\mu\notin\mathbb{B},\; \xi\in H^C,\; \big(\mu I-\cP^{*C}\big)\xi=0,\; \big(D(\cd)\bigl|_{Z}\big)^{*C}\Psi_a(T,\cd)^{*C}\xi=0\;\;\mbox{over}\;\; (0,T)\Rightarrow \xi=0.
\end{equation}
Suppose that $\mu$ and $\xi$ satisfy the conditions on the left side of (\ref{a.2.1}).
Write $\xi=\xi_1+i\xi_2$ where $\xi_1, \xi_2\in H$. Then, we have
\begin{equation}\label{YUHUANG5.10}
\big(D(\cd)\bigl|_{Z}\big)^*\Psi_a(T,\cd)^{*}\xi_j=0,\; j=1,2.
\end{equation}
Clearly,  $\psi_j(\cdot)\triangleq\Psi_a(T,\cd)^{*}\xi_j$ (with $j=1,2$) is the solution to the  equation (\ref{a.2.7}) where $n_0T$ and $\xi$ are replaced by $T$ and $\xi_j$ respectively.
Since $\Psi_a(T,\cd)^*\xi_j$ is continuous on $[0,T]$ and $D(t)$ is independent of $t$, it follows from (\ref{YUHUANG5.10}) that
\begin{equation}\label{a.2.2}
\big(D(0)\bigl|_{Z}\big)^{*}\Psi_a(T,0)^{*}\xi_j=0,\; j=1,2.
\end{equation}
For each $\eta\in H$, we have $P^*\eta\in P^*H$. This, along with  (\ref{a.2.3}), yields
$$\ba{rl}
\ns &\left\langle \big(D(0)\bigl|_{Z}\big)^{*}\Psi_a(T,0)^{*}\xi_j,~P^{*}\eta\right\rangle=\left\langle \Psi_a(T,0)^{*}\xi_j,~\big(D(0)\bigl|_{Z}\big)P^{*}\eta\right\rangle\\
\ns&%
=\left\langle \Psi_a(T,0)^{*}\xi_j,~\chi_{\omega}P^{*}\eta\right\rangle=\left\langle P\chi_{\omega}\Psi_a(T,0)^{*}\xi_j,~\eta\right\rangle,\; j=1,2.
\ea
$$
This, alone with (\ref{a.2.2}), implies that
$P\chi_{\omega}\Psi_a(T,0)^{*}\xi_j=0,\; j=1,2$, from which, we have
\begin{equation}\label{a.2.4}
\lan P^{*}\Psi_a(T,0)^{*}\xi_j,\,\chi_{\omega}\Psi_a(T,0)^{*}\xi_j\ran=\lan \Psi_a(T,0)^{*}\xi_j,\,P\chi_{\omega}\Psi_a(T,0)^{*}\xi_j\ran=0,\;j=1,2.
\end{equation}
Tow facts are as follows. First, it follows from (\ref{2-4}) that
\begin{equation}\label{a.2.5}
 P^{*}\Psi_a(T,0)^{*}\xi_j=\Psi_a(T,0)^{*}P^{*}\xi_j,\;j=1,2.
 \end{equation}
Second, by (\ref{WANGJIAXU3.31}), (\ref{a.1.50}), and the first three conditions on the left side of (\ref{a.2.1}), we have $\xi\in \tilde H_1^C$.
Since $P^*=\tilde P$ and $\tilde P$ is a projection from $H$ to $\tilde H_1$ (see Lemma~\ref{HUANGlemma2.3}), we see that
$P^*: H\rightarrow\tilde H_1$ is a projection. Hence, $P^{*C}: H^C\rightarrow \tilde H_1^C$ is a projection. These two facts yields that  $ P^{*C}\xi=\xi$, from which, it follows that
$P^*\xi_j=\xi_j$, $j=1,2$. This
along with  (\ref{a.2.4}) and (\ref{a.2.5}), indicates that
$\left\|\chi_{\omega}\Psi_a(T,0)^{*C}\xi\right\|=0$, i.e., $\chi_\omega\psi_j(T)=0$.
By the unique continuation  property of parabolic equations established in \cite{LFH} (see also \cite{Phung1} and \cite{Phung2}), we find that $\xi_j=0$, $j=1,2$, which leads to $\xi=0$.
   This completes the proof.
\endpf

 \vskip 10pt

Finally, we will present a controlled heat equation which is not LPFS. Write
$\l_1$ and $\l_2$ for the first and the second eigenvalues of  the operator $-\triangle$ with $\cD(-\triangle)=H_0^1(\Omega)\cap H^2(\Omega)$, respectively. Let $\xi_j$, $j=1,2$,  be an eigenfunction corresponding to $\l_j$.
Consider the following heat equation:
\begin{equation}\label{h1state} \left\{\begin{array}{lll}
\ns \partial_ty(x,t)-\triangle y(x,t)-\l_2y(x,t)=\lan u(t),\xi_1\ran\xi_1(x) &\mbox{in}&Q,\\
\ns y(x,t)=0&\mbox{on}&\sum.
\end{array}\right.\end{equation}
where $u(\cdot)\in L^2(\mathbb{R}^+;L^2(\Omega))$.
By a direct calculation, one has that
$$
V_{n_0}=\span\{\xi_1\}\;\;\mbox{and}\;\; H_1\supseteq\mbox{span}\{\xi_1,\xi_2\}.
$$
These, along with the equivalence $(a)\Leftrightarrow (b)$ in Theorem \ref{theorem1}, indicates that  (\ref{h1state}) is not LPFS.

\section*{Appendix}
\setcounter{equation}{0}
 \renewcommand{\theequation}{
            a.\arabic{equation}}

\noindent{\it The proof of Lemma \ref{lemma2.1}.}
By the  compactness of  $\{S(t)\}_{t>0}$, the assumption $(\cH_2)$ and (\ref{huang1.4}), one can easily check that  each  $\cP(t)$, with $t\geq 0$,
is  compact.  Hence,  each
 $\cP(t)^C: H^C\rightarrow H^C$, with $t\geq 0$, is also compact. Thus,  for each $t\geq 0$, $\sigma(\cP(t)^C)\setminus\{0\}$ consists
of all nonzero eigenvalues  $\{\l_j(t)\}^\infty_{j=1}$ (in $\mathbb{C}$) of $\cP(t)^C$ such that
$\lim\limits_{j\rightarrow\infty}|\l_j(t)|=0$.

We next show that $\{\l_j(t)\}^\infty_{j=1}$ is independent of $t$.
For this purpose, we arbitrarily fix
 $s_1$ and $s_2$ with  $0\le s_1\leq s_2+ T$. Let $\l\in \mathbb{C}^1$ be a non-zero eigenfunction of $\cP(s_1)^C$ and $\eta\in H^C$ be a corresponding eigenfunction, i.e.,
\begin{equation}\label{wgs2.1}
\cP(s_1)^C\eta=\l\eta.
\end{equation}
 Write  $\l=\a_1+i\a_2$ with $\a_1, a_2\in \mathbb{R}$ and $\eta=\eta_1+i\eta_2$ with $\eta_1, \eta_2\in H$. By (\ref{wgs2.1}), we have
\begin{equation}\label{wgs2.2}
\cP(s_1)\eta_1=\a_1\eta_1-\a_2\eta_2,\qq\cP(s_1)\eta_2=\a_2\eta_1+\a_1\eta_2.
\end{equation}
From  (\ref{wgs2.2}) and (\ref{wgs1.5}), one can easily check that
$\cP(s_2)^C(\Phi(s_2,s_1)^C\eta)=\l\Phi(s_2,s_1)^C\eta$.
This implies that $\l$ is an eigenvalue of $\cP(s_2)^C$ and $\Phi(s_2,s_1)^C\eta$ is a corresponding eigenfunction. Hence,
$$\sigma(\cP(s_1)^C)\setminus\{0\}\subseteq\sigma(\cP(s_2)^C)\setminus\{0\}.$$
Similarly, we can show
$$\sigma(\cP(s_2)^C)\setminus\{0\}\subseteq\sigma(\cP(s_1+T)^C)\setminus\{0\}.$$
Then by the $T$-periodicity of $\cP(\cdot)$, $\sigma(\cP(t)^C)\setminus\{0\}$ is independent of $t$.
This completes the proof.
 \endpf

\noindent{\it The proof of Lemma \ref{lemma2.2}.} First of all, we let
\begin{equation}\label{wgs2.7}
\hat H_1(t)\triangleq \hat P(t)H^C\;\;\mbox{and}\;\; \hat H_2(t)\triangleq (I-\hat P(t))H^C,\;\; t\geq 0.
\end{equation}
From Theorem 6.17 on Page 178 in \cite{Kato}, it follows that when $t\geq 0$,
both $\hat H_1(t)$ and $\hat H_2(t)$ are invariant w.r.t. $\cP(t)^C$;
\begin{equation}\label{YUHUANG2.7}
\hat P(t): H^C(t)\rightarrow \hat H_1(t)\;\;\mbox{is a projection};
\end{equation}
\begin{equation}\label{2-6}H^C= \hat H_1 (t)\bigoplus \hat H_2 (t);
\end{equation}
 and
 \begin{equation}\label{2-7}
\sigma\big(\cP(t)^C|_{\hat H_1(t)}\big) =\{\l_j \}^n_{j=1}\;\;\mbox{and}\;\; \sigma(\cP(t)^C|_{\hat H_2 (t)})\setminus\{0\}= \{\l_j \}^\infty_{j=n+1},
\end{equation}
where $\{\l_j\}_{j=1}^\infty$ and $n$ are given by (\ref{WGS1.8}) and (\ref{wgs1.9}) respectively.

Then we prove that the operator $P(t)$, with $t\geq 0$, is a linear operator from $H$ to $H$. For this purpose,  it suffices to show that
\begin{equation}\label{wgs2.11}
\hat P(t)h\in H,\;\;\mbox{when}\;\; h\in H\;\;\mbox{and}\;\; t\geq 0.
\end{equation}
The proof of (\ref{wgs2.11}) is as follows. By (\ref{wgs1.11}),  it holds that
\begin{equation}\label{wgs2.12}
 \hat P(t)h=\ds \frac{-\d}{2\pi}\int^{2\pi}_{0} \big( \d e^{i\theta} I-\cP(t)^C\big)^{-1} e^{i\theta}d\theta\ h,\,\;\;\mbox{when}\;\; h\in H\;\;\mbox{and}\;\; t\geq 0.
 \end{equation}
Write $F$ for  the conjugate map from $H^C$ to $H^C$, i.e., $F(h+ig)=h-ig$ for any $h,g\in H$.
We claim
\begin{equation}\label{wgs2.13}
F\big(\big(\delta e^{i\theta}I-\cP(t)^C\big)^{-1}e^{i\theta} h\big)=(\delta e^{-i\theta}I-\cP(t)^C)^{-1}e^{-i\theta}h\;\;\mbox{for all}\;\;\theta\in [0,2\pi], h\in H\;\;\mbox{and}\;\; t\geq 0,
\end{equation}
When (\ref{wgs2.13}) is proved, it follows from (\ref{wgs2.12}) and (\ref{wgs2.13}) that
 $$\ba{rl}
 F(\hat P(t))h &=\ds  \frac{-\d}{2\pi}\int^{2\pi}_{0} \left( \d e^{-i\theta} I-\cP(t)^C\right)^{-1}e^{-i\theta}d\theta\,h\\
  &=\ds \frac{-\d}{2\pi}\int^{2\pi}_{0} \left(\d e^{i\theta} I-\cP(t)^C\right)^{-1}
 e^{i\theta}d\theta\,h=\hat P(t)h\;\;\mbox{for each}\;\; t\geq 0, h\in H,
 \ea$$
 Which leads to (\ref{wgs2.11}). Now we are on the position to show (\ref{wgs2.13}). Arbitrarily fix $\theta\in [0, \pi]$, $t\geq 0$ and $h\in H$. Write
 \begin{equation}\label{wgs2.14}
 (\delta e^{i\theta}I-\cP(t)^C)^{-1}e^{i\theta}h=g_1+ig_2,\;\; g_1,g_2\in H.
 \end{equation}
 It is clear that $(\delta e^{i\theta}I-\cP(t)^C)(g_1+ig_2)=e^{i\theta}h$, from which,
 one can directly check that
$$
(\delta e^{-i\theta}I-\cP(t)^C)(g_1-ig_2)=e^{-i\theta}h.
$$
Hence, $(\delta e^{-i\theta}I-\cP(t)^C)^{-1}(e^{-i\theta}h)=g_1-ig_2=F(g_1+ig_2)$.
This, along with (\ref{wgs2.14}), leads to (\ref{wgs2.13}).

Next we prove that $P(t)$, with $t\geq 0$, is  a projection on  $H$.
Let $H_1(t)$ and $H_2(t)$, with $t\geq 0$, be defined by (\ref{2-2-1}). Two observations are given in order:
\begin{equation}\label{2-8}
 \ba{rl}\ns\hat H_1(t)&\triangleq\hat P(t)H^C=\big\{\hat P(t)(h_1+ih_2)\bigm|h_1,h_2\in H\big\}=\big\{ P(t)h_1+iP(t)h_2\bigm|h_1,h_2\in H\big\}\\
 \ns&=P(t)H+i P(t)H\triangleq H_1(t)+iH_1(t)\triangleq H_1^C(t); \ea
 \end{equation}
  \begin{equation}\label{2-9}
 \hat H_2(t)= H_2^C(t).
 \end{equation}
By (\ref{YUHUANG2.7}) and (\ref{2-8}), we see that  $\hat P(t)$ (with $t\geq 0$) is a projection from $H^C$ onto $H_1(t)^C$. Thus, for each  $t\geq 0$,
$$
P(t)(h_1+h_2)= \hat P(t)(h_1+h_2)=h_1,\;\;\mbox{when}\;\; h_1\in H_1(t), h_2\in H_2(t),
$$
i.e.,  $P(t)$ is a projection from $H$ onto $H_1(t)$.
Besides,  (\ref{2-2}) follows from (\ref{2-6}), (\ref{2-8}) and (\ref{2-9}).

Finally,  we will  show properties $(a)$-$(f)$ one by one.

\noindent {\it The proof of $(a)$:}
Since $\cP(\cdot)$ is $T$-periodic,  so is $\hat P(\cdot)$ (see (\ref{wgs1.11})). This, along with (\ref{P}), indicates the  $T$-periodicity of  $P(\cdot)$. Then by (\ref{2-2-1}), both $H_1(\cdot)$ and $H_2(\cdot)$ are $T$-periodic.

\noindent {\it The proof of $(b)$:}  Let $t\geq 0$. Since $\hat H_1(t)$ and $\hat H_2(t)$ are invariant w.r.t. $\cP(t)^C$,
so are $H_1(t)^C$
and $H_2(t)^C$ (see (\ref{2-8}) and (\ref{2-9})).
 Hence, $H_1(t)$ and $H_2(t)$ are invariant w.r.t. $\cP(t)$.

\noindent {\it The proof of $(c)$:} (\ref{{2-3}}) follows from (\ref{2-7}), (\ref{2-8}) and (\ref{2-9}). Meanwhile, by (\ref{number}) and (\ref{{2-3}}), we see that
$\mbox{dim} H_1(t)^C=n_0$,
 which leads
to  (\ref{wgs1.14}).

\noindent {\it The proof of $(d)$ and $(e)$:} Let $0\leq s\leq t<\infty$. By (\ref{wgs1.5}), we have that $\Phi(t,s)\cP(s)=\cP(t)\Phi(t,s)$. From this, one can directly verify that
$\Phi(t, s)^C\hat P(s)=\hat P(t)\Phi(t, s)^C$.
 This, along with (\ref{2-2-1}), (\ref{P}) and (\ref{wgs2.11}), indicates that
\begin{equation}\label{2-11}
\Phi(t, s)H_1(s)\subseteq P(t)H\triangleq H_1(t),
\end{equation}
which leads to $(e)$. Meanwhile, it follows from (\ref{2-11}) that $\Phi(t,s)\in \mathcal{L}(H_1(s), H_1(t))$. Similarly, one can show that $\Phi(t,s)\in \mathcal{L}(H_2(s), H_2(t))$. Hence, $(d)$ stands.

\noindent{\it The proof of $(f)$:} Let $\bar \rho\triangleq (-{\ln \bar\d})/{T}>0$ with $\bar\delta$ given by (\ref{decay}). Because of (\ref{{2-3}}), it follows from Theorem 4 on Page 212 in \cite{Yosida}
 that
the spectral radius of $\cP^C(0)\bigl|_{H_2(0)^C}$  equals to $\bar\delta$. Thus, we have
$$\bar\d=\lim\limits_{k\rightarrow\infty}\big\|\big(\cP(0)^C\bigl|_{H_2(0)^C}\big)^k\big\|^{\frac{1}{k}}.$$
Now we arbitrarily fix a  $\rho\in(0,\bar\rho)$ where $\bar \rho$ is given by (\ref{decay}). Then it holds that
$\bar\d\triangleq e^{-\bar\rho T}<e^{-\rho T}$.
Thus there is positive integer $\hat N$ such that
$\big\|\big(\cP(0)^C\bigl|_{H_2(0)^C}\big)^k\big\|<e^{-\rho kT}$ for all $k\ge \hat N$,
which implies
\begin{equation}\label{2-14-1}
\big\|\big(\cP(0)\bigl|_{H_2(0)}\big)^k\big\|<e^{-\rho kT}\;\; \mbox{ for all }\;\;k\ge \hat N.
\end{equation}

Notice that $\Phi(\cdot, \cdot)$ is continuous from $[0,T]\times[0,T]$ to $\cL(H)$ (see Lemma 5.6 on Page 68 in \cite{Li}). Thus, we can write
\begin{equation}\label{wgs2.19}
C_1\triangleq\max\limits_{0\le t_1\le t_2\le T}\left\|\Phi(t_2,t_1)\right\|\in \mathbb{R}^+; \;\; C_\rho\triangleq(C_1+1)^2e^{3\rho T}\in \mathbb{R}^+.
\end{equation}
We are going to show that the above $C_\rho$ satisfies (\ref{2-5}). For this purpose, we let $0\leq s\leq t<\infty$ and $h_2\in H_2(s)$. For each $r\in \mathbb{R}^+$, we denote by
$[r]$ the integer such that $r-1<[r]\leq r$. There are only two possibilities: $(i)$ $[t/T]=[s/T]$ and $(ii)$ $[t/T]\neq [s/T]$.
In the first  case, it follows from (\ref{wgs1.5}) and (\ref{wgs2.19}) that
$$\ba{rl}&\left\|\Phi(t,s)h_2\right\|=\left\|\Phi(t-\left[{s}/{T}\right]T,s-\left[{s}/{T}\right]T)h_2\right\|
\le \left\|\Phi(t-\left[{s}/{T}\right]T,s-\left[{s}/{T}\right]T)\right\|\|h_2\|\\
\ns&\le C_1\|h_2\|\le C_1e^{\rho T}e^{-\rho(t-s)}\|h_2\|<(C_1+1)^2e^{3\rho T}e^{-\rho(t-s)}\|h_2\|=C_\rho e^{-\rho(t-s)}\|h_2\|,
\ea$$
i.e., $C_\rho$ satisfies  (\ref{2-5}) in the first  case. In the second case, we have that  $[t/T]T\geq [s/T]T+T$; and  it follows from  $(d)$ and $(a)$ that
$\Phi([s/T]T+T,s)h_2\in H_2(0)\triangleq H_2$. These, along with  (\ref{wgs1.5}) and (\ref{2-14-1}), indicate that
$$\ba{rl}
&\left\|\Phi(t,s)h_2\right\|\leq \big\|\Phi\big(t,\,\big[{t}/{T}\big]T\big)\big\|\cdot
\big\|\cP(0)^{([{t}/{T}]-[{s}/{T}]-1)}\Phi\big(\big[{s}/{T}\big]T+T,\,s\big)h_2\big\|\\
\ns&\le\big\|\Phi\big(t-\big[{t}/{T}\big]T,\,0\big)\big\|\cdot e^{-\rho T([{t}/{T}]-[{s}/{T}]-1)}\cdot\big\|
\Phi\big(T,\,s-\big[{s}/{T}\big]T\big)\big\|\|h_2\|.
\ea$$
By this and (\ref{wgs2.19}), one can directly check that
$\big\|\Phi(t,s)h_2\big\|\le  C_\rho e^{-\rho (t-s)}\|h_2\|$, i.e.,  $C_\rho$ satisfies (\ref{2-5}) in the second case. This shows (\ref{2-5}) and completes the proof.
  \endpf

\bigskip

\noindent {\it The proof of Lemma \ref{HUANGlemma2.3}.}
By (\ref{a.1.41}), (\ref{a.1.43}), (\ref{a.1.42}), (\ref{YUHEHHEN2.22}) and (\ref{YUHEHHEN2.23}),  one can make use of the exactly same way utilized in the proof of Lemma~\ref{lemma2.2} to verify all properties in Lemma~\ref{HUANGlemma2.3}, except for (\ref{a.1.49})-(\ref{WANGJIAXU3.31}).
Since $(\l I-\cP^{*C})^{-1}=\big((\l I-\cP^C)^{-1}\big)^*$, (\ref{a.1.49}) follows from (\ref{a.1.43}),  (\ref{wgs1.11}) and (\ref{P}).
Now, we  prove (\ref{a.1.50}).
The first equation of (\ref{a.1.50}) follows from the definition of $\tilde H_1$ and (\ref{a.1.49}).
It is clear that $P^*H\supseteq P^*H_1$.
On the other hand, since
$P^*P h=0\Rightarrow\lan h, P^*P h\ran=0\Rightarrow P h=0$, it holds that $\cN(P^*P)\subseteq \cN(P)$. This, together with the fact that $H_1=PH$ (see (\ref{2-2}) and (\ref{wgs1.15})), yields

$$ P^*H_1=P^*PH=\cR(P^*P)=\cN(P^*P)^\bot\supseteq\cN(P)^\bot=\cR(P^*)=P^*H.$$
Therefore, (\ref{a.1.50}) holds.

The proof of (\ref{WANGJIAXU3.31}) is as follows.
Since $\cP^{*C}\xi=\mu\xi$, we derive from (\ref{a.1.43}) that
$$
\hat {\tilde P}\xi=\ds\frac{1}{2\pi i}\int_{\Gamma}\left(\l I-\cP^{*C}\right)^{-1}d\l\,\xi=
\ds\frac{1}{2\pi i}\int_{\Gamma}\left(\l -\mu\right)^{-1}d\l\,\xi=\xi.
$$
Hence, $\xi\in \hat {\tilde P}H^C$. Meanwhile, by the definitions of $\tilde P$ and $\tilde H_1$,  we find that
$$
\tilde H_1^C=\big( \hat {\tilde P}\big|_H\big)^CH^C= \hat {\tilde P} H^C.
$$
Thus, it holds that $\xi\in\tilde H_1^C$. This completes the proof.

\endpf

\vskip 5pt

\noindent {\bf Acknowledgement} The authors would like to  thank Professor Xu Zhang for his valuable suggestions on this paper.

\end{document}